\documentclass[a4paper,11pt]{article}
\pdfoutput=1
\usepackage{mathrsfs}
\usepackage[utf8]{inputenc}
\usepackage[T1]{fontenc}
\usepackage{amsmath}
\usepackage{amssymb}
\usepackage{amsthm}
\usepackage[francais,english]{babel}
\usepackage{graphicx}
\usepackage{color}
\usepackage{fullpage}
\usepackage{dsfont}
\usepackage{enumerate}
\usepackage{caption}

\definecolor{purple}{rgb}{0.57,0.1,0.53}

\newenvironment{hyp}[1]{\par\bigskip \noindent \textbf{Assumption #1 : }\itshape} {\par\bigskip }
\newcommand{\HypLevy}{A}
\newcommand{\HypMut}{B.2}
\newcommand{\HypMutConst}{B.1}

\newtheorem{lemme}{Lemma}[section]
\newtheorem{theoreme}[lemme]{Theorem}

\newtheorem{proposition}[lemme]{Proposition}

\newtheorem{remarque}[lemme]{Remark}
\newenvironment{demo}{\noindent\emph{Proof :} \\ }{\hfill $\square$ \\ \par\smallskip}
\newenvironment{demoth}[1]{\noindent\emph{Proof of Theorem #1 :} \\ }{\hfill $\square$}
\newenvironment{demopr}[1]{\noindent\emph{Proof of Proposition #1 :} \\ }{\hfill $\square$}

\newenvironment{demole}[1]{\noindent\emph{Proof of Lemma #1 :} \\ }{\hfill $\square$}
\newcounter{claim}

\newcommand{\Z}{\mathbb Z}
\renewcommand{\L}{\Lambda}
\newcommand{\tL}{\tilde\Lambda}

\newcommand{\sachant}{\,|\,}
\newcommand{\R}{\mathbb R}
\newcommand{\N}{\mathbb N}
\newcommand{\F}{\mathcal F}
\renewcommand{\d}[1]{\text{d}#1}
\renewcommand{\P}{\mathbb P}
\newcommand{\E}{\mathbb E}
\newcommand{\point}{\,\cdot\,}
\newcommand{\eps}{\varepsilon}

\newcommand{\D}{\mathbb D}

\newcommand{\Exc}{\mathscr E}

\newcommand{\Linv}{L^{-1}}
\newcommand{\Tx}{{T^{-x}}}
\newcommand{\Txn}{{T_{n}^{-x}}}
\renewcommand{\kill}{k}

\newcommand{\Zbn}{\tilde Z^{\textsc{m}}_n} 
\newcommand{\Zb}{Z^{\textsc{m}}}

\newcommand{\Zc}{\hat Z^{\textsc{m}}} 
\newcommand{\Zcn}{\hat Z_n^{\textsc{m}}}

\newcommand{\Han}{H^{+}_n}

\newcommand{\Hcn}{H^{\textsc{m}}_n}
\newcommand{\mathcalHcn}{\mathcal H^{\textsc{m}}_n}
\newcommand{\Ha}{H^{+}}

\newcommand{\Hc}{H^{\textsc{m}}}

\newcommand{\teta}{\tilde\eta}

\newcommand{\Fc}{\mathcal F^{\textsc{m}}}

\newcommand{\B}[1]{\mathbb B_{#1}}

\renewcommand{\ll}{\textit{\texttt l}}

\newcommand{\en}{\text{e}_n}
\newcommand{\e}{\text{e}}

\newcommand{\zero}{\{0\}}
\newcommand{\un}{\{1\}}

\newcommand{\Tree}{\mathbb T}
\newcommand{\tTree}{\tilde{\mathbb T}}

\newcommand{\Mponct}{\mathcal M_\textsc{P}}
\newcommand{\tZ}{\tilde Z}
\newcommand{\tZntau}{\tilde Z_{n,\tau}}
\newcommand{\exc}{e}

\newcommand{\Leb}{\text{Leb}}

\renewcommand{\ni}{{}_n^{(i)}}
\newcommand{\neps}{{}_{n,\eps}}
\newcommand{\nieps}{{}_{n,\eps}^{(i)}}

\newcommand{\Ups}{\Upsilon}
\newcommand{\MUps}{\Upsilon^\textsc{m}}
\newcommand{\MUpsn}{\Upsilon_n^\textsc{m}}

\newcommand{\m}{\textsc{m}}

\newcommand{\num}{\nu^\textsc{m}}
\newcommand{\nud}{\nu^{\textsc{d}}}
\newcommand{\nui}{\nu_\eps^{\textsc {init}}}
\newcommand{\nuin}{\nu_{n,\eps}^{\textsc {init}}}

\newcommand{\sigmad}{\sigma^{\textsc{k}}}
\newcommand{\Hd}{H^{\textsc{k}}}
\newcommand{\mud}{\mu^{\textsc{k}}}
\newcommand{\nuk}{\nu^{\textsc{k}}}
\newenvironment*{remerciements}{%

\begin{abstract}
}{\end{abstract}}

\begin{document}

\selectlanguage{english}
\title{Lévy processes with marked jumps II : \\ Application to a population model with mutations at birth}
\author{Cécile Delaporte\footnote{UPMC Univ. Paris 6, Laboratoire de probabilités et modèles aléatoires CNRS UMR 7599, and Collège de France, Center for interdisciplinary research in biology CNRS UMR 7241 ; \textit{e-mail :} cecile.delaporte@upmc.fr}}
\date{}
\maketitle

\begin{abstract}
Consider a population where individuals give birth at constant rate during their lifetimes to i.i.d. copies of themselves. Individuals bear clonally inherited types, but (neutral) mutations may happen at the birth events. The smallest subtree containing the genealogy of all the extant individuals at a fixed time $\tau$ is called the coalescent point process. We enrich this process with the history of the mutations that appeared over time, and call it the marked coalescent point process.

With the help of limit theorems for Lévy processes with marked jumps established in \cite{LPWM1}, we prove the convergence of the marked coalescent point process with large population size and two possible regimes for the mutations - one of them being a classical rare mutation regime, towards a multivariate Poisson point process.  This Poisson point process can be described as the coalescent point process of the limiting population at $\tau$, with mutations arising as inhomogeneous regenerative sets along the lineages. Its intensity measure is further characterized thanks to the excursion theory for spectrally positive Lévy processes. In the rare mutations asymptotic, mutations arise as the image of a Poisson process by the ladder height process of a Lévy process with infinite variation, and in the particular case of the critical branching process with exponential lifetimes, the limiting object is the Poisson point process of the depths of excursions of the Brownian motion, with Poissonian 
mutations on the lineages.
\end{abstract}


\selectlanguage{french}
\section{Introduction}\label{sec_Intro}
\qquad A splitting tree (\cite{Geiger}, \cite{GK}, \cite{ALContour}) describes a population of individuals with i.i.d. lifetime durations, which distribution is not necessarily exponential, giving birth at constant rate during their lives. Each birth gives rise to a single child, who behaves as an independent copy of her parent. We consider here the extended framework of \cite{ALContour} : for each individual, the birth times and lifetimes of her progeny is given by a Poisson process with intensity $\d t\cdot\L(\d r)$, where the so-called lifespan measure $\L$ is a Lévy measure on $(0,\infty)$ satisfying $\int(1\wedge r)\L(\d r)<\infty$. In particular, the number of children of a given individual is possibly infinite. In addition, we assume that individuals carry types, and that every time a birth occurs, a mutation may happen, giving rise to a mutant child. Mutations are assumed to be neutral, meaning that they do not affect the behaviour of individuals. In order to take this into account, we introduce 
\textit{marked splitting trees} : to each birth event we associate a mark in $\{0,1\}$, which will code for the absence ($0$) or presence ($1$) of a mutation. In other words, a $0$-type birth means a clonal birth, and a $1$-type birth produces a mutant child. The mutations experienced by the population are then described by these marks.

Population models with mutations have inspired lots of works in the past, and have many applications in domains such as population genetics, phylogeny or epidemiology. Such models have been well studied in the particular case of populations with fixed size. In the Wright-Fisher and Moran models with neutral mutations, as well as in the Kingman coalescent, explicit results on the allelic partition of the population are provided by Ewens' sampling formula (\cite{Ewens},\cite{Durrett}). Relaxing the hypotheses of constant population size, branching processes with mutations at birth are studied in the monography \cite{Taib}. More recently, results have been obtained for the allelic partition and frequency spectrum of splitting trees, with mutations appearing either at birth of individuals (\cite{Richard}) or at constant rate along the lineages (\cite{LambertSpecies}, \cite{CL1}, \cite{CL2}), and are reviewed in \cite{CLR}. 
The present work focuses on asymptotic results when the size of the population gets large, for the genealogy (with mutational history) of splitting trees with mutations at birth, and relies on a previous article \cite{LPWM1}.

\paragraph*{Genealogy of the $n$-th population} \quad 
Let us fix some positive real number $\tau$. For $n\in\N$, consider a marked splitting tree $\Tree_n$, and condition it on having a fixed positive number $I_n$ of individuals alive at level $\tau$. Note that we use here the word 'level' to denote the real time in which the individuals live, whereas we reserve the word 'time' for the index of stochastic processes. 
This paper follows on from a work of L. Popovic (\cite{Popovic}) in the critical case with exponential lifetimes, without mutations, in which she proved the convergence in distribution of the coalescent point process (i.e. the smallest subtree containing the genealogy of the extant individuals) towards a certain Poisson point process. Our aim is to provide asymptotic results as $I_n$ gets large, for the structure of the genealogy of the population up to level $\tau$, enriched with the random levels at which marks occurred on the lineages. 
To this aim, after a proper rescaling of $\Tree_n$, we introduce a random point measure $\Sigma_n$ which we call the marked coalescent point process. This point measure has $I_n-1$ atoms ; its $i$-th one is itself a random point measure, whose set of atoms contains all the levels where mutations occurred on the $i$-th lineage, and the coalescence time between individuals $i$ and $i-1$. 
This sequence of point measures $(\Sigma_n)$ is the mathematical object for which we aim to get convergence as $n\to\infty$, after having set some convergence assumptions, which we discuss later.

Our work mainly relies on the study of splitting trees with the help of the so-called jumping chronological contour process (or JCCP). This process is an exploration process of the tree (without mutations) introduced by A. Lambert in \cite{ALContour}, visiting all the existence levels of all the individuals exactly once, and ending at level $0$. 
He showed in this paper that the JCCP of a tree truncated up to level $\tau$ is a compensated compound Poisson process with no negative jumps (spectrally positive Lévy process with finite variation) reflected below $\tau$ and killed when hitting $0$.  In particular, the labeling of the excursions of the JCCP below $\tau$ provides a labeling of the extant individuals at level $\tau$. Inferring properties concerning the genealogy of the alive population at level $\tau$ in the tree then essentially consists in studying the excursions away from $\tau$ of this reflected Lévy process. 

We introduce in \cite{LPWM1} a generalization of this contour process to the framework of our rescaled marked splitting trees $(\tTree_n)$. We are thereby led to study a bivariate Lévy process $(\tZ_n,\tZ_n^\m)$. Roughly speaking, $\tZ_n$ codes for the JCCP of $\tTree_n$ (without mutations), and $\tZ_n^\m$ codes for the mutations. Namely, since a jump of $\tZ_n$ corresponds to the encounter of a birth event when exploring the tree, $\tZ_n^\m$ will jump as well (with amplitude 1) if this birth was of type $1$. The process $(\tZ_n,\tZ_n^\m)$ is in one-to-one correspondence with the marked tree $\tTree_n$. We now want to characterize the law of the atoms of $\Sigma_n$ using this property. 
Let us first give an idea of our reasoning in the case where there is no mutations. The JCCP of $\tTree_n$, truncated up to $\tau$, is distributed as $\tZ_n$ reflected below $\tau$. The set of levels at which births occurred on the lineage of the $i$-th individual, up to its coalescence with the rest of the tree, is then exactly the set of values taken by the future infimum of the $i$-th excursion of the JCCP under $\tau$. 
First, this entails that the atoms of $\Sigma_n$ are i.i.d. Second, using a time reversal argument, the distribution of this set can be read from the ascending ladder height process of $\tZ_n$. 
A similar reasoning for the splitting tree with mutations leads to the following facts. Consider $\Ha_n$ the ascending ladder height process of $\tZ_n$, and put marks on its jumps in agreement with the marks on the corresponding jumps of $\tZ_n$. Note that this implies a selection of the marks that are carried by jumps of the supremum process of $\tZ_n$. 
Denoting by $\Hc_n$ the counting process of these marks, the bivariate process $(\Ha_n,\Hc_n)$ is a (possibly killed) bivariate subordinator which we call the marked ladder height process. The mutations on a lineage form then an inhomogeneous regenerative set, distributed as the image by $\Ha_n$ of the jump times of $\Hc_n$ under the excursion measure of $\tZ_n$ away from $0$, which finally yields a simple description of the law of the (i.i.d.) atoms of $\Sigma_n$.

\paragraph*{Convergence results} \quad 
Obtaining an invariance principle for a population model in a large population asymptotic requires to assume that as $n\to\infty$, the population converges in a certain sense.  A classical example would be the convergence of the rescaled Bienaymé-Galton-Watson process towards the Feller diffusion (\cite{Lamperti}). Now regardless of mutations, the JCCP offers a one-to-one correspondence between $\Tree_n$ and a continuous time process. Our first assumption arises then naturally as the convergence in distribution of the properly rescaled Lévy process $\tZ_n$ towards a Lévy process $Z$ (with infinite variation, Assumption \HypLevy). In particular, the lifetimes of individuals do not necessarily vanish in the limit. 
Besides, two different assumptions concerning the mutations are considered. The first one (\HypMutConst) falls within the classical asymptotic of rare mutations : every birth is of type 
$1$ with a constant probability $\theta_n$, and $\theta_n\to0$ as $n\to\infty$. Asymptotic results in this framework are obtained in \cite{Bertoin2} for the genealogical structure of alleles in a critical or subcritical Bienaymé-Galton-Watson process (however contrary to ours, they do not concern the extant population at a fixed time horizon, but the whole population). The second one (\HypMut) examines the case where the probability of an individual to be a mutant is correlated with her lifetime, in the sense that mutations favor longer lifetimes. 

While Assumption \HypLevy\ alone ensures the convergence in distribution of $\Ha_n$ towards the classical ladder height process of $Z$, Assumptions \HypMutConst\ and \HypMut\ are designed to allow that of the marked ladder height process. Indeed, we prove in \cite{LPWM1} the convergence in law of $(\Ha_n,\Hcn)$ towards a (possibly killed) bivariate subordinator $(\Ha,\Hc)$, such that $\Ha$ is the ladder height process of $Z$. Note nevertheless that in this framework there is in general no convergence of the whole mutation process, namely $\tZ_n^\m$. 
In the case of Assumption \HypMutConst, $\Ha$ and $\Hc$ are independent, and $\Hc$ is a Poisson process with parameter $\theta$, which arises as the limit of the sequence $\theta_n$ after a proper rescaling. This means that the contribution to the mutations in the limit exclusively comes from individuals with vanishing lifetimes. This is no longer the case under Assumption \HypMut, yet additional independent marks can appear if $Z$ has a Gaussian component. 
Using this convergence to deduce that of the (rescaled) law of the mutations on a lineage, the convergence of $(\Sigma_n)$ to a Poisson point measure is then a straightforward consequence of the law of rare events for null arrays (see e.g. \cite[Th. 16.18]{Kall}). Under \HypMutConst, its intensity measure is the law of the image by $\Ha$ of an independent Poisson process with parameter $\theta$, under the excursion measure of $Z$ away from zero. A very similar but slightly more complicated result, involving the limiting marked ladder height process $(\Ha,\Hc)$, is available under \HypMut. Besides, in the case where $Z$ is a Brownian motion, $\Ha$ is simply a drift, and thus the intensity measure is the law of a Poisson process killed at some independent random time, distributed as the depth of an excursion of the Brownian motion away from 0.

\paragraph{Outline} \quad 
The paper is organized as follows : Section \ref{sec_Prelim} sets up notation for the topological framework, and provides some background on the excursion theory for Lévy processes (see e.g. \cite{B} and \cite{K}). Section \ref{sec_Results} is devoted to the statement of our results, and Section \ref{sec_Proof} to their proofs. In the appendix, we give proof of some properties that are consequences of Assumption \HypLevy, and which we make frequent use of throughout the paper.

\setlength{\parindent}{0cm}
\section{Preliminaries}\label{sec_Prelim}
\subsection{Topology} \label{sec_topo}

We consider the Euclidean space $\R^d$ and endow it with its Borel $\sigma$-field $\mathcal B(\R^d)$. We denote by 
$\D(\R^d)$ the space of all càd-làg functions from $\R_+$ to $\R^d$. We endow the latter with the Skorokhod topology, which makes it a Polish space (see \cite[VI.1.b]{JS}). In the sequel, for any function $f\in\D(\R)$ and $x>0$, we will use the notation $\Delta f(x)=f(x)-f(x-)$, where $f(x-)=\lim_{u\to x,\,u<x} f(u)$.\\

Now for any Polish space $X$, with its Borel $\sigma$-field $\mathcal B$, the space $\mathcal M_f(X)$ of positive finite measures on $(X,\mathcal B)$ can be endowed with the weak topology :
It is the coarsest topology for which the mappings $\mu\mapsto\int f\d\mu$ are continuous for any continuous bounded function $f$. In the sequel, we will use the notation $\mu(f):=\int f\d\mu$.\\

Hence we endow here $\mathcal M_f(\R^d)$ and $\mathcal M_f(\D(\R^d))$ with their respective weak topologies. The notation $\Rightarrow$ will be used for both weak convergence in $\R^d$ and in $\D(\R^d)$, and we will use the symbol $\stackrel{(d)}{=}$ for the equality in distribution. Recall that for any sequence of $\R^d$-valued càd-làg processes $(X_n)$, the weak convergence of $(X_n)$ towards a process $X$ of $\D(\R^d)$ is equivalent to the finite dimensional convergence of $(X_n)$ towards $X$ along any dense subset $D\subset\R_+$, together with the tightness of $(X_n)$. For more details about convergence in distribution in $\D(\R^d)$, see \cite[VI.3]{JS}.\\

From now on, we fix $\tau>0$. We consider the space of positive point measures on $(0,\tau)\times\{0,1\}$, and endow it with the $\sigma$-field generated by the mappings $\{p_B:\xi\mapsto \xi(B),\ B\in\mathcal B((0,\tau))\otimes\mathscr P(\{0,1\})\}$. Then we denote by $\Mponct$ the subset of the point measures on $(0,\tau)\times\{0,1\}$ of the form 
$$\delta_{(a_m,0)}+\sum_{i=0}^{m-1}\delta_{(a_i,1)}, \text{ where }m\in\Z_+ \text{ and } 0<a_0<\ldots<a_{m-1}\leq a_m<\tau.$$ 
The trace $\sigma$-field on $\Mponct$ is in particular generated by the class 
$$\mathscr C=\big\{p_B,\ B=[\eps,\tau)\times\zero,\ B=[\eps,\tau)\times\un\big\}_{\eps\in(0,\tau)}.$$

\subsection{Excursion theory for spectrally positive Lévy processes} \label{sec_SPLP}

We provide here some background about the excursion theory for spectrally positive Lévy processes. For the basic properties concerning spectrally positive Lévy processes that will be needed here, we refer to a summary we provide in \cite[Section 2]{LPWM1} (these properties can otherwise be found in \cite{B} or \cite{K}, for example). \\

Let $X$ be a spectrally positive Lévy process with Lévy measure $\L$. We define its past supremum $\bar X_t:=\underset{[0,t]}\sup X$ for all $t\geq0$. We denote by $\Exc$ the set of excursions of $X-\bar X$ away from $0$ : $\Exc$ is the set of the càd-làg functions $\epsilon$ with no negative jumps for which there exists $\zeta=\zeta(\epsilon)\in(0,\infty]$, which will be called the lifetime of the excursion, and such that $\epsilon(0)=0$, $\epsilon(t)$ has values in $(-\infty,0)$ for $t\in(0,\zeta)$ and in the case where $\zeta<\infty$, $\epsilon(\zeta)\in[0,\infty)$.\\

The reflected process $X-\bar X$ is a Markov process for which one can construct $(L_t)_{t\geq0}$ a local time at 0. We denote by $\Linv$ its inverse, and we consider the process $e=(e_t)_{t\geq0}$ with values in $\Exc\cup \{\partial\}$ (where $\partial$ is an additional isolated point), defined by :
\[e_t:=\left\{ 
\begin{array}{l l}
  ((X-\bar X)_{s+\Linv(t-)},0\leq s<\Linv(t)-\Linv(t-)) & \quad \text{if } \Linv(t-)<\Linv(t)\\
  \partial & \quad \text{else}\\ \end{array} \right. \]

Then according to Theorem IV.10 in \cite{B}, if $X$ does not drift to $-\infty$, then $0$ is recurrent for the reflected process, and $(t,e_t)_{t\geq0}$ is a Poisson point process with intensity $c\;\d t\ N(\d\epsilon)$, where $c$ is some constant depending on the choice of $L$, and $N$ is a measure on $\Exc$. Else, $(t,e_t)_{t\geq0}$ is a Poisson point process with intensity $c\;\d t\ N(\d\epsilon)$, stopped at the first excursion with infinite lifetime. \\

In the same way, we denote by $\Exc'$ the set of excursions of $X$ away from $0$ : $\Exc'$ is the set of the càd-làg functions $\epsilon$ with no negative jumps for which there exists $\zeta=\zeta(\epsilon)\in(0,\infty]$, and such that $\epsilon(0)=0$, $\epsilon(t)$ has values in $\R^*$ for $t\in(0,\zeta)$, and $\epsilon(\zeta)=0$ if $\zeta<\infty$. We then introduce $\chi(\epsilon):=\inf\{t\in(0,\zeta],\ \epsilon(t)\in[0,\infty)\}$. 

Denoting by $\mathscr L$ a local time at $0$ of $X$ and by $\mathscr L^{-1}$ its inverse, we define the process $e'=(e'_t)_{t\geq0}$ with values in $\Exc'\cup \{\partial\}$
\[e'_t:=\left\{ 
\begin{array}{l l}
  (X_{s+\mathscr L^{-1}(t-)},0\leq s<\mathscr L^{-1}(t)-\mathscr L^{-1}(t-)) & \quad \text{if } \mathscr L^{-1}(t-)<\mathscr L^{-1}(t)\\
  \partial & \quad \text{else}\\ \end{array} \right. \]
If $X$ has no Gaussian component, any excursion $e'_t\in\Exc'$ first visits $(-\infty,0)$, and we necessarily have $\chi(e'_t)>0$ (but possibly infinite). On the other hand, if $X$ has a Gaussian component, it can creep upwards and then $\chi(e'_t)\in[0,\infty]$.
 Again, according to Theorem IV.10 in \cite{B}, $e'$ is a Poisson point process with intensity $c'\;\d t\ N'(\d\epsilon)$, stopped if $X$ is subcritical at the first excursion with infinite lifetime. Here $c'$ is some constant depending on the choice of $\mathscr L$ and $N'$ a measure on $\Exc'$.\\

Finally, we describe some marginals of $N$ and $N'$ in the proposition below, for which we refer to \cite[Th. 6.15 and (8.29)]{K}, \cite[(3)]{BertoinDecomp} and \cite[Cor. 1]{BertoinPitman}.\\

\begin{proposition} \label{prop_Exc_mesure} 
 We have for all $z,x>0$ :
\begin{enumerate}[\upshape(i)]
 \item If $X$ has finite variation,
\upshape$$N(-\epsilon(\zeta-)\in \d x,\ \epsilon(\zeta)\in \d z,\ \zeta<\infty)=W(0)e^{-\eta x}\d x\L(x+\d z)$$
\itshape 
 \item  If $X$ has infinite variation and no Gaussian component (i.e. $b=0$),
\upshape$$N(-\epsilon(\zeta-)\in \d x,\ \epsilon(\zeta)\in \d z,\ \zeta<\infty)=e^{-\eta x}\d x\L(x+\d z).$$
\end{enumerate}
\itshape Moreover, in both cases, under $N(\,\cdot\,\sachant-\epsilon(\zeta-)=x,\ \zeta<\infty)$, the reversed excursion 
$$\big(-\epsilon((\zeta-t)-),\ 0\leq t<\zeta\big)$$ 
is equal in law to $(X_t,\ 0\leq t<T^0)$ under $\P_x(\,\cdot\,\sachant T^0<\infty)$. \\

Finally, the same statement holds replacing $N$ by $N'$ and $\zeta$ by $\chi$. 
\end{proposition}

\section{A limit theorem for splitting trees with mutations at birth}
\label{sec_Results}

\subsection{JCCP of a marked splitting tree} \label{Sec_ST_with_mutations}

Formally, a splitting tree (without mutations) is a random real tree characterized by a $\sigma$-finite measure $\Lambda$ on $(0,\infty)$, satisfying $\int(1\wedge u) \L(\d u)<\infty$. Consider such a splitting tree, and assume first that there is extinction of the population. In \cite{ALContour}, A. Lambert considers a contour process of this tree called JCCP (jumping chronological contour process). He establishes that the tree and its contour process are in one-to-one correspondence and characterizes the law of the latter : conditional on the first individual in the tree to have life duration $x$, its JCCP is distributed as a finite variation, spectrally positive Lévy process with drift $-1$ and Lévy measure $\L$, starting at $x$, and killed upon hitting $0$. In the case of non extinction, we then can consider the JCCP of the tree truncated up to level $\tau$, which has the law of the Lévy process described above, starting at $x\wedge\tau$, and reflected below level $\tau$. 
As noticed in Section \ref{sec_Intro}, the exploration of the tree by its JCCP defines a way of ordering the individuals. In the sequel, when we label the extant individuals at level $\tau$, we refer to that order.\\

Consider now a marked splitting tree $\Tree$ as defined in Section \ref{sec_Intro}. We assume that the probability for a child to be a mutant can only (possibly) depend on her life span $u$, and if we denote by $f(u)$ this probability, where $f$ is a function from $\R_+^*$ to $[0,1]$, $f$ will be called the\textit{ mutation function} of the tree. Then $\Tree$ is characterized by its mutation function $f$ and its lifespan measure $\L$.\\

\begin{figure}[!h]
\begin{center}
  \includegraphics[width=15.5cm]{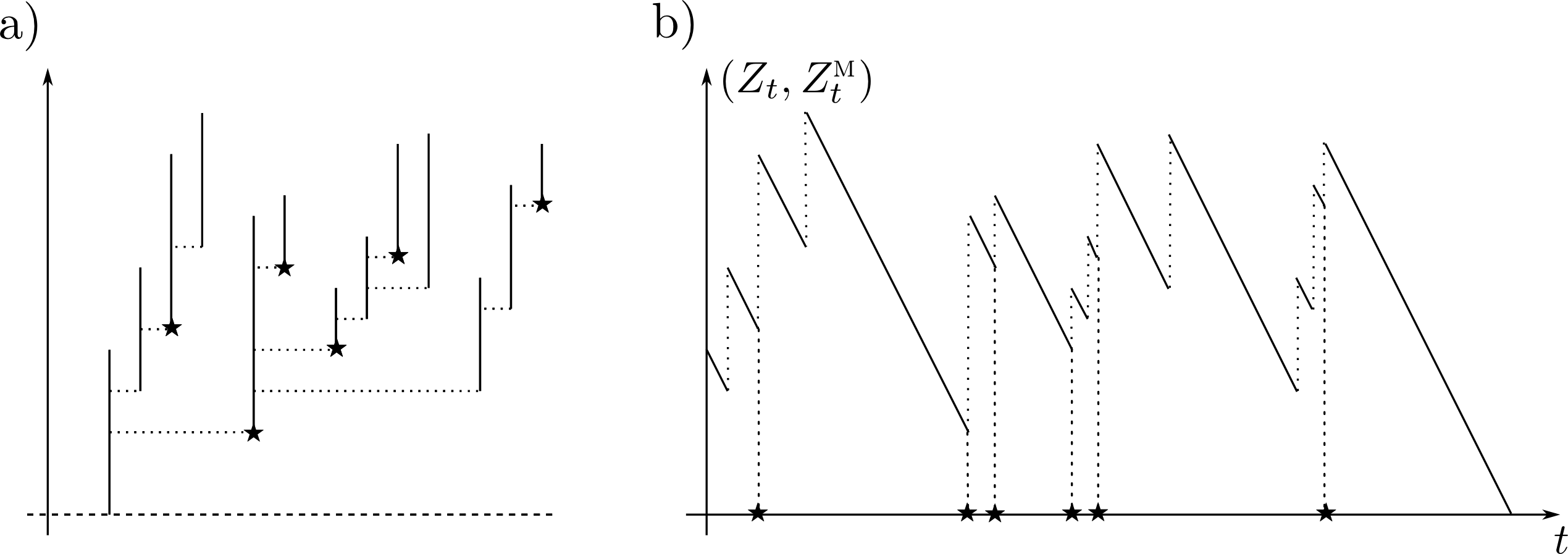}
  \caption{a) A marked splitting tree : the vertical axis indicates chronological levels ; the horizontal axis has no meaning, but the horizontal lines show filiation. The marks in the tree are symbolized by stars.\\
b) The associated JCCP $(Z,Z^{\textsc m})$ : $Z$ is the classical JCCP of the (non-marked) tree represented in a) ; the counting process of the mutations $Z^{\textsc m}$ is not drawn as a jump process on $\R_+$, but is represented by the sequence of its jump times, which are symbolized by stars on the horizontal axis.} 
\label{figure_jccp}
\end{center}
 \end{figure}

Then similarly as in the case without mutations, we define the JCCP of $\Tree$. First assume that there is extinction of its population. Then the JCCP of the marked tree $\Tree$ is a bivariate process $(Z,Z^{\textsc m})$ from $\R_+$ to $\R_+\times\Z_+$, whose first coordinate $Z$ is the JCCP of the splitting tree without marks, and whose second coordinate $Z^{\textsc m}$ is the counting process of the mutations (see Figure \ref{figure_jccp}). 
More precisely, for every jump time of $Z$ (which corresponds to the encounter of a birth event in the exploration process), $Z^{\textsc m}$ jumps (with amplitude $1$) iff this birth was a $1$-type birth. Hereafter we say that a jump of $Z$ occurring at time $t$ carries a mark (or a mutation) if 
$\Delta Z^\m (t)=1$.\\

This bivariate process is in one-to-one correspondence with $\Tree$. Besides, conditional on the first individual to have life duration $x$, it is distributed as a bivariate Lévy process with drift $(-1,0)$, and Lévy measure $\L(\d u)\B{f(u)}(\d q)$ (where $\B r$ denotes the Bernoulli probability measure with parameter $r$), starting at $(x,0)$, and killed as soon as its first 
coordinate hits $0$.  As in the non-marked case, if the assumption of extinction does not hold, the law of the JCCP of the truncated tree can be obtained from the Lévy process we just described.\\

\subsection{Definitions and notation}\label{sec_Defs}

\subsubsection{Rescaling the population}
Let $(\L_n)_{n\geq1}$ be a sequence of measures on $(\R_+^*,\mathcal B(\R_+^*))$ satisfying $\int(1\wedge u) \L_n(\d u)<\infty$ for all $n$, and $(f_n)_{n\geq1}$ a sequence of continuous functions from $\R^+$ to $[0,1]$. \\

We now consider a sequence of marked splitting trees $(\Tree_n)_{n\geq1}$ such that for all $n$, $\Tree_n$ has lifespan measure $\L_n$, and mutation function $f_n$. 
Recalling that  $\B r$ denotes the Bernoulli probability measure with parameter $r$, we consider $(Z_n,Z_n^{\textsc m})$ an independent bivariate Lévy process with finite variation, Lévy measure $\L_n(\d u)\B{f_n(u)}(\d q)$ and drift $(-1,0)$, and make the following assumption : 

 \begin{hyp}\HypLevy
  There exists a sequence of positive real numbers $(d_n)_{n\geq1}$ such that as $n\to\infty$, the process defined by
$$\tilde Z_n:=\Big(\frac1n Z_{n}(d_n t)\Big)_{t\geq0}$$ 
converges in distribution to a (necessarily spectrally positive) Lévy process $Z$ with infinite variation. We denote by $\L$ its Lévy measure and by $b$ its Gaussian coefficient ($b\in\R_+$).
 \end{hyp}

For all $n\in\N$ and for all $t\geq0$, set $\tilde Z_n^{\textsc m}(t):=Z_n^{\textsc m}(d_n t)$. With an abuse of notation, the law of $(\tilde Z_n,\tilde Z_n^{\textsc m})$ conditional on $(\tilde Z_n(0),\tilde Z_n^\m(0))=(x,0)$, and the law of $Z$ conditional on $Z(0)=x$, will both be denoted by $\P_x$, and we write $\P$ for $\P_0$.\\ 

Denote by $\tTree_n$ the splitting tree obtained from $\Tree_n$ by rescaling the branch lengths by a factor $\frac1n$. The introduction of the process $(\tilde Z_n,\tilde Z_n^\m)$ is motivated by its fundamental role in the characterization of the law of the JCCP of $\tTree_n$ truncated up to level $\tau$ (see later Lemma \ref{lemme_exc_indep}).

\paragraph*{Some notation :} 
The Laplace exponents $\psi_n$ of $Z_n$, $\tilde \psi_n$ of $\tilde Z_n$ and $\psi$ of $Z$ are defined by 
$$\E(e^{-\lambda Z_n(t)})=e^{t\psi_n(\lambda)},\  \E(e^{-\lambda \tilde Z_n(t)})=e^{t\tilde\psi_n(\lambda)}\  \text{ and }\ \E(e^{-\lambda Z(t)})=e^{t\psi(\lambda)},\ \ \ \lambda\geq0.$$
We denote by $\teta_n$ (resp. $\eta$) the largest root of $\tilde\psi_n$ (resp. $\psi$) and by $\tilde\phi_n$ (resp. $\phi$) the inverse of $\tilde\psi_n$ (resp. $\psi$) on $[\teta_n,\infty)$ (resp. $[\eta,\infty)$). We denote by $\tilde W_n$ (resp. $W$) the scale function of $\tilde Z_n$ (resp. $Z$). Finally, we denote by $\tL_n$ the Lévy measure of $\tilde Z_n$.

\paragraph*{Remarks about $(d_n)$ :} 
Writing for $\lambda\geq0$, $\E(e^{-\lambda \tilde Z_n(t)})=e^{d_n t\psi_n(\lambda/n)}$, we get from the Lévy-Khintchine formula \cite[(2)]{LPWM1} that $\tilde Z_n$ has drift $-\frac{d_n}n$, Lévy measure $\tL_n=d_n\L_n(n\cdot)$ and Laplace exponent $\tilde\psi_n=d_n\psi(\cdot/n)$. In particular, this gives $\tilde W_n(0)=n/d_n$. We prove in the appendix that $\tilde W_n$ converges pointwise to $W$ as $n\to\infty$, and besides, the assumption of infinite variation of $Z$ ensures $W(0)=0$. Thereby we know that necessarily $\frac{d_n}{n}\to\infty$ as $n\to\infty$.

\subsubsection{Asymptotic for the mutations}

In order to allow the convergence in distribution of the mutation levels on the lineages, we have to make some technical assumptions on the mutation functions $f_n$. Here we suggest two possible assumptions : in the first one, the probability of a child in $\Tree_n$ to be a mutant is constant, while in the second one, this probability depends on its life duration.

\begin{hyp}\HypMutConst
\begin{enumerate}[\upshape(a)]
 \item For all $n\geq1$, for all $u\in\R_+$, $f_n(u)=\theta_n$, where $\theta_n\in[0,1]$.
 \item As $n\to\infty$, $\frac{d_n}{n}\theta_n$ converges to some finite real number $\theta$.
 \end{enumerate}
\end{hyp}

\begin{hyp}\HypMut
\begin{enumerate}[\upshape(a)]
 \item The sequence $\big(u\mapsto \frac{f_n(nu)}{1\wedge u}\big)$ converges uniformly to $u\mapsto \frac{f(u)}{1\wedge u}$ on $\R_+^*$. \label{HypMut_cv}
 \item There exists $\kappa\geq0$ such that $f(u)/u \to \kappa$ as $u\to0^+$. \label{HypMut_0} 
\end{enumerate}
\end{hyp}

Note that in \HypMutConst, necessarily $\theta_n\to 0$ as $n\to\infty$, corresponding to the classical rare mutation asymptotic. Then if we denote by $f$ the limit of the sequence $(f_n)$, we have $f\equiv0$. Besides, in Assumption \HypMut\ the choice of $f_n$ and $f$ is independent of $\tilde Z_n$ and $Z$. \\

\begin{remarque}\label{remark_pas_de_cv_Zm}
These two possible assumptions for the rescaling of the mutations have been chosen so that as $n\to\infty$, the marked coalescent point process converges. However this choice does not imply, despite Assumption \HypLevy, the convergence of the bivariate process $(\tilde Z_n,\tilde Z_n^\m)$. As pointed out in \cite{LPWM1}, it is even never the case under \HypMut.
\end{remarque}

\subsubsection{Marked genealogical process}

From now on, we consider the sequence of rescaled marked splitting trees $(\tTree_n)$, and condition $\tTree_n$ on having $I_n$ extant individuals at level $\tau$, where $I_n\sim\frac{d_n}{n}$ as $n\to\infty$.\\

Consider a realization of $\tTree_n$, and label the $I_n$ individual alive at $\tau$ from $0$ to $I_n-1$ (according to Section \ref{Sec_ST_with_mutations}). Then to the $i$-th one we associate a simple point measure $\sigma\ni$, with values in $(0,\tau)\times\{0,1\}$, as follows :\\
Consider the lineage of individual $i$, and assume it contains $M$ $1$-type birth events. Denote by $m_0$ the level where the lineage coalesces with the rest of the tree, and by $m_j,\ 1\leq j\leq M$ the successive levels (in increasing order) where the 1-type birth events happened. Then we set 

$$\sigma\ni:=\delta_{(\tau-m_0,0)}+\sum_{1\leq j\leq M} \delta_{(\tau-m_j,1)}.$$

\begin{figure}[!h]
\begin{center}
  \includegraphics[width=15.5cm]{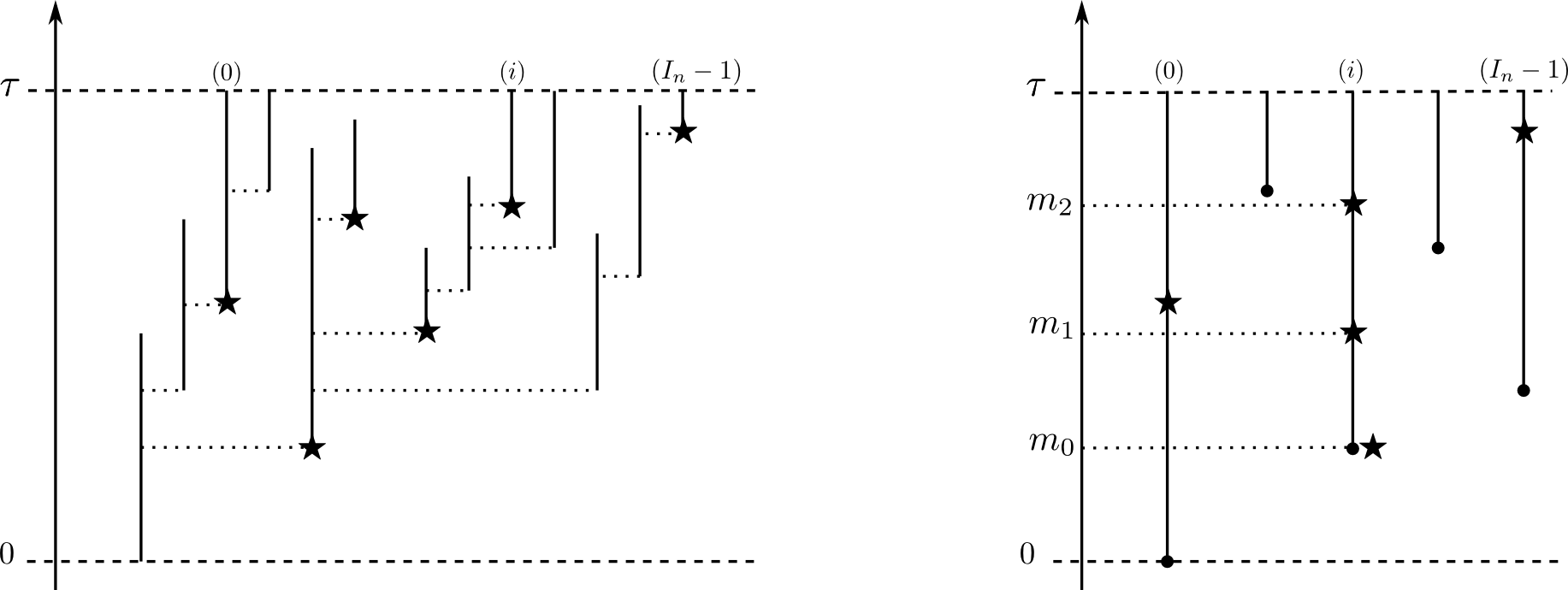}
 \caption{A marked splitting tree truncated up to level $\tau$ and the associated marked coalescent point process. The $1$-type birth events are symbolized by stars, and dots represent coalescence levels. In this example, the coalescence between the lineages of individuals $i$ and $i-1$ coincides with a $1$-type birth event, and we have $\sigma\ni=\delta_{(\tau-m_0,0)}+\delta_{(\tau-m_0,1)}+\delta_{(\tau-m_1,1)}+\delta_{(\tau-m_2,1)}.$}
\label{figure_sigmani}
\end{center}
 \end{figure}

Hence the point measure $\sigma\ni$ is in the space $\Mponct$, and keeps record of all the mutation events on the $i$-th lineage, and of the coalescence level of this lineage with the rest of the tree (see Figure \ref{figure_sigmani}). The quantity $\tau-m_0$ will be called the coalescence time of the lineage (the word 'time' is here to interpret as a duration). Note that in case the coalescence corresponds to a 1-type birth event, we have $m_0=m_1$.\\

Now for all $n\geq1$, we define the following random point measure on $[0,1]\times\Mponct$ :
$$\Sigma_n :=\sum_{i=1}^{I_n-1} \delta_{\{\frac{in}{d_n},\sigma\ni\}}.$$

The first individual (labeled $0$) is on purpose not taken in account (see Remark \ref{remark_first_lineage} below). The point measure $\Sigma_n$ is called the marked coalescent point process of $\tTree_n$. As announced in Section \ref{sec_Intro}, the aim of this paper is to obtain a convergence theorem for $\Sigma_n$ in a large population asymptotic. 

\subsection{Main results}\label{sec_Main_results}

\begin{figure}[!h]
\begin{center}
  \includegraphics[width=15.5cm]{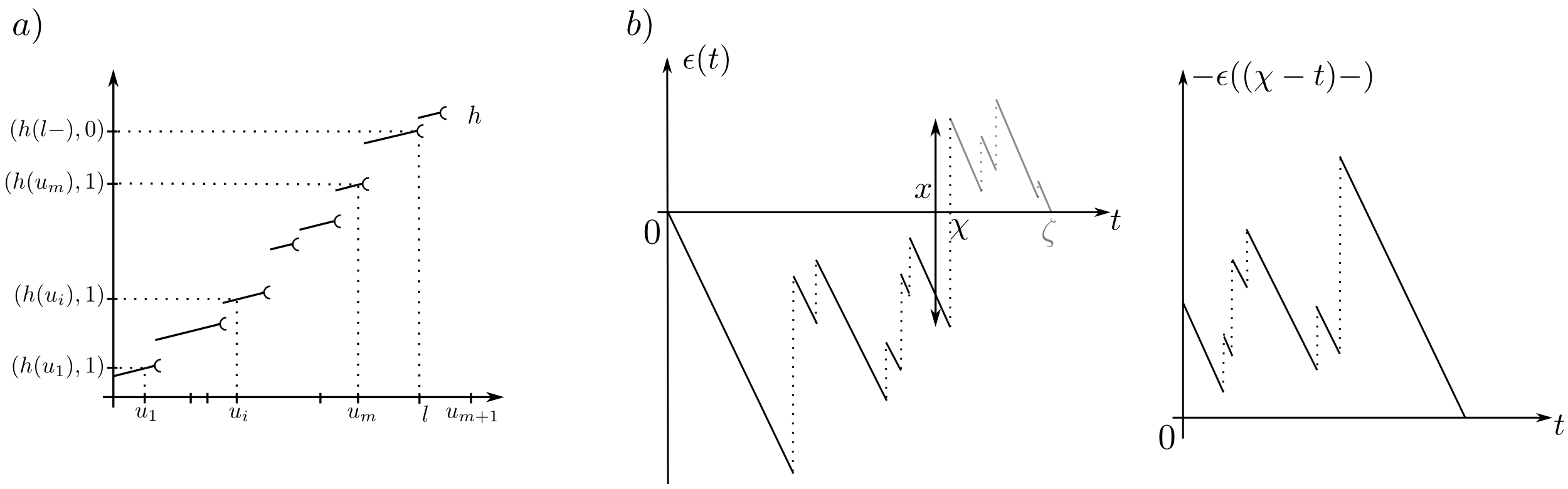}
 \caption{
a) A graphical representation of a triplet $(h,\textbf{u}=(u_i)_{i\geq1},l)\in\D((0,\tau))\times(\R_+)^{\N}\times\R_+$. In this example, $j(\textbf{u},l)=m$ and $\Psi(h,\textbf{u},l)=\delta_{(h(l-),0)}+\sum_{i=1}^{m} \delta_{(h(u_i),1)}$. \\
b) Left panel : A representation (in finite variation) of an excursion $\epsilon\in\mathscr E'$, with finite lifetime $\zeta$, such that $\Delta\epsilon(\chi)=x$. Right panel : The corresponding reversed excursion on $[0,\chi)$ : $(-\epsilon(\chi-t)-))_{0\leq t<\chi}$ (which belongs to $\mathscr E''$).}
\label{figure_N}
\end{center}
 \end{figure}

We first introduce some notation. To begin with, we define the mapping $\Psi$ as follows (see figure \ref{figure_N}.a) : for all  $(h,\textbf{u}=(u_i)_{i\geq1},l)\in\D((0,\tau))\times(\R_+)^{\N}\times\R_+$, 

$$\Psi(h,\textbf{u},l) = \delta_{(h(l-),0)}+\sum_{i=1}^{j(\textbf{u},l)} \delta_{(h(u_i),1)},
$$
where
\[j:\left\{\begin{array}{ccc}
   (\R_+)^{\N}\times\R_+ &\to &\N\cup\{+\infty\} \\
   (\textbf{u},l) &\mapsto &\sup\{i\geq1,\ u_i\leq l\}
  \end{array}\right.\]

The function $\Psi$ has values in the point measures on $(0,\tau)\times\{0,1\}$, and if $j(\textbf{u},l)<+\infty$ and $h(u_1)<\ldots<h(u_{j(\textbf{u},l)})\leq h(l-)$, then $\Psi(h,\textbf{u},l)$ is in the set $\Mponct$.\\

For any càd-làg piecewise-constant function $g:\R_+\to(0,\tau)$, if $(g_i)_{i\geq1}$ denotes the sequence of its jump times (with $g_1=0$ in case $g(0)>0$), we will use the notation $\Psi(h,g,l)$ instead of $\Psi(h,(g_i),l)$. \\

We denote by $\bar Z(t):=\sup_{[0,t]} Z$ the current supremum process of $Z$, and by $\Ha:=\bar Z\circ \Linv$ the ladder height process of $Z$, where $L$ is a local time at the supremum for $Z$, which will be specified later (see Section \ref{sec_cvMLHP}), and $\Linv$ its inverse local time. We denote by $T^A$ the first entrance time of $Z$ in the Borel set $A$, and write  $T^x$ for $T^{\{x\}}$. \\

Finally, we denote by $N'$ the excursion measure of $Z$ away from zero (see Section \ref{sec_SPLP}), and we choose the normalization of the local time $\mathscr L$ according to \cite{Obloj}, i.e. $\mathscr L$ satisfies the equality $\E\big(\int_{(0,\infty)}e^{-t}\d\mathscr L_t\big)=\phi'(1)$. Recall that for $\epsilon\in\mathscr E'$, $\chi(\epsilon)$ denotes its first entrance time into $[0,\infty)$. Define $\mathscr E''$ the set of all càd-làg functions $\epsilon$ with lifetime $\zeta<\infty$, such that $\epsilon(0)=\epsilon(\zeta)=0$ and $\epsilon(x)>0$ for all $0<x<\zeta$. Then we define a measure $N''$ on $\mathscr E''\times\{0,1\}$ as follows (see Figure \ref{figure_N}.b) : for all $(\epsilon,\epsilon^\m)\in\mathscr E''\times\{0,1\}$, 
$$N''((\epsilon,\epsilon^\m)\in\,\d E\times \d q):=\int_{[0,\infty)} N'(\Delta\epsilon(\chi)\in\d x,\ (-\epsilon(\chi-t)-)_{0\leq t<\chi}\in\d E)\B{f(x)}(\d q).$$
Note that 
$$N''((\epsilon,\epsilon^\m)\in\,\d E\times \{0,1\})=N'(\,(-\epsilon(\chi-t)-)_{0\leq t<\chi}\in\d E),$$
and that in the case where $Z$ does not drift to $+\infty$, the excursions of $Z$ have finite lifetime, and from a time reversal argument we have for any measurable set $B$ of $\mathscr E''$,  $N''(B\times\{0,1\})=N'(\epsilon|_{[\chi,\zeta)}\in B)$, where $\epsilon|_{[\chi,\zeta)}$ denotes the restriction of $\epsilon$ to the interval $[\chi,\zeta)$.

\paragraph*{Results under Assumption \HypMutConst}\quad

In this paragraph we suppose that Assumptions \HypLevy\ and \HypMutConst\ are satisfied.

\begin{theoreme}\label{th_cv_PPP_Const}
Consider an independent Poisson process $\Theta$ with parameter $\theta$. We introduce $\sigma$, a random element of $\Mponct$, defined on $\{T^0<\infty\}$ by 
 $$\sigma=\Psi(\Ha,\Theta,L(T^0)).$$
Then the sequence $(\Sigma_n)$ converges in distribution towards a Poisson point measure $\Sigma$ on $[0,1]\times\Mponct$ with intensity measure 
\upshape $\Leb\otimes\Pi_1,$ \itshape
where  \upshape Leb \itshape denotes the Lebesgue measure, and $\Pi_1$ is a measure on $\Mponct$ defined by
$$\Pi_1=N''(\sigma\in\,\cdot\,,\ \sup\epsilon<\tau).$$
\end{theoreme}

\begin{remarque}\label{remark_mesure_infinie}                                                                                 
Denote by $B_{\geq m}:=\{\sigma\in\Mponct,\ \sigma((0,\tau)\times\un)\geq m\}$ the set of point measures of $\Mponct$ having at least $m$ points with second coordinate $1$ in the interval $(0,\tau)$, which can be interpreted here as the presence of at least $m$ mutations on a lineage. Then the measure $\Pi_1(B_{\geq1})$ is not necessarily finite (see Example 1).
\end{remarque}

\begin{remarque} \label{remark_first_lineage}
 Note that we excluded in $\Sigma_n$ the first lineage $\sigma_n^{(0)}$, for which without additional assumption, we cannot easily get a similar result as for the other lineages. However, if we assume that the lifetime of the first individual in $\tTree_n$ converges as $n\to\infty$ towards some value greater than $\tau$, we can adapt Theorem \ref{th_cv_PPP_Const}. The limiting object is then obtained by adding to $\Sigma$ a Dirac mass on $(0,\delta_{(\tau,0)})$.
\end{remarque}

\begin{remarque} \label{remark_conditioning}\textbf{Conditioning $\tTree_n$ on survival at level $\tau$} \\
We obtain a similar result if, instead of conditioning $\tTree_n$ on having $I_n$ extant individuals at level $\tau$, we condition it on survival at level $\tau$. Indeed, if we denote by $\tilde\Xi_n(\tau)$ the number of extant individuals in $\tTree_n$ at level $\tau$, we know that conditional on $\tilde\Xi_n(\tau)\geq1$, $\tilde\Xi_n(\tau)$ follows a geometric distribution with parameter $\frac{n}{d_n\tilde W_n(\tau)}$ (see \cite[prop.5.6]{ALContour} ). Then thanks to the pointwise convergence of $\tilde W_n$ towards $W$ (see Proposition \ref{prop_cv}), we get that $\frac{n}{d_n}\tilde\Xi_n(\tau)$ converges in distribution towards an exponential variable with parameter $\frac1{W(\tau)}$. \\
Then the sequence $(\Sigma_n)$ converges in law to a Poisson point measure on $[0,\textbf{e}]\times\Mponct$ with intensity \upshape Leb\itshape$\otimes\Pi_1$, where $\textbf{e}$ is an independent exponential variable with parameter $\frac1{W(\tau)}$.
\end{remarque}

Assume now that $Z$ has \textbf{no Gaussian component}, and let $\Theta$ be as in Theorem \ref{th_cv_PPP_Const}. Then using Proposition \ref{prop_Exc_mesure} we get
\[\Pi_1=\int_{(0,\tau)} \d x\,\bar\L(x)\,\P_x(\sigma\in\,\cdot\,,\ T^0<T^{(\tau,\infty)} ),
\]
where $\sigma$ is  defined in Theorem \ref{th_cv_PPP_Const} and $\bar\L(x)=\L((x,\infty))$ for all $x>0$. Hence in the limit, the mutations appearing on a lineage are distributed according to a point measure $\sigma$, where $\sigma$ is distributed as the image of the jump times of an independent Poisson process with parameter $\theta$, by the ladder height process of $Z$ conditioned on $T^0<T^{(\tau,\infty)}$, and starting at the opposite of the undershoot of an excursion with depth smaller than $\tau$. \\

Finally, the following proposition expresses the law of $\sigma$ under $\P_x(\,\cdot\,\cap\, \{T^0<T^{(\tau,\infty)}\})$ in terms of the image of an independent Poisson process by an inhomogeneous killed subordinator.
\begin{proposition}\label{prop_mesure_de_sauts}
Let $\Theta$ and $\sigma$ be as in Theorem \ref{th_cv_PPP_Const}. For all $x\in(0,\tau)$, 
$$\P_x(\sigma\in\,\cdot\,,\ T^0<T^{(\tau,\infty)})=\P_x(\sigmad\in \,\cdot\,),$$
where 
$$\sigmad:=\Psi(\Hd,\Theta,L),$$
with $\Hd$ a killed inhomogeneous subordinator with drift $\frac{b^2}{2}$ and jump measure $\mud$, defined for all $a\in(0,\tau)$ and $u\in(0,\tau-a)\times\{+\infty\}$ by :
\upshape$$\mud(a,\d u):=\frac{1}{W(a)}\delta_{+\infty}(\d u)+ \int_{(0,a)} \d x\,\L(x+\d u)\,\frac{W(a-x)W(\tau-a-u)}{W(a)W(\tau)},$$\itshape
and $L:=\inf\{t\geq0,\ \Hd(t)=+\infty\}$ the killing time of $\Hd$.
\end{proposition}

\paragraph*{Results under Assumption \HypMut}\quad

We suppose now that Assumptions \HypLevy\ and \HypMut\ are satisfied. We establish in this case some very similar results as under \HypMutConst, but in a slightly more complicated version. Indeed, Assumption \HypMutConst\ ensures the independence of $\Ha$ with a certain process we define later (namely the subordinator $\Hc$ that appears in the following statement), while in case \HypMut\ these two subordinators are no longer independent.

\begin{theoreme}\label{th_cv_PPP}
There exists a process $\Hc$, starting at $0$ under $N''$, such that $(\Ha,\Hc)$ is a (possibly killed) bivariate subordinator, and such that $(\Sigma_n)$ converges in distribution towards a Poisson point measure $\Sigma$ on $[0,1]\times\Mponct$ with intensity measure 
\upshape $\Leb\otimes\Pi_2,$ \itshape
where 
$$\Pi_2=N''(\,\Psi(\Ha,\epsilon^\m+\Hc,L(T^0))\,\in\,\cdot\,,\ \sup\epsilon<\tau).$$
The processes $\Ha$ and $\Hc$ are not independent unless $Z$ is a Brownian motion with drift, and the law of $(\Ha,\Hc)$ is explicitly characterized in Theorem \ref{th_cv_H}.
\end{theoreme}

Note that Remarks \ref{remark_first_lineage} and \ref{remark_conditioning} are still relevant in case \HypMut.

\begin{remarque} \label{remark_Brownien_B2} 
If the limiting process $Z$ is a Brownian motion with drift, $\Ha$ is a deterministic drift and hence $\Ha$ and $\Hc$ are automatically independent. Hence in this case, Theorem \ref{th_cv_PPP_Const} remains valid under Assumption \HypMut.
\end{remarque}

Similarly as under \HypMutConst, if $Z$ has no Gaussian component we can reexpress the measure $\Pi_2$ as follows :
\[\Pi_2=\int_{(0,\tau)\times\{0,1\}} \d x\,\int_{(x,\infty)} \L(\d u)\,\B{f(u)}(\d q)\,\P_x(\sigma_q\in\,\cdot\,,\ T^0<T^{(\tau,\infty)} ),
\]
where for $q\in\{0,1\}$, $\sigma_q=\Psi(\Ha,q+\Hc,L(T_0))$. \\

Furthermore, as in Proposition \ref{prop_mesure_de_sauts}, we have for all $x\in(0,\tau)$, $q\in\{0,1\}$
$$\P_x(\sigma_q\in\,\cdot\,,\ T^0<T^{(\tau,\infty)} )=\P_x(\sigmad_q\in\,\cdot\,),$$
where
$$\sigmad_q:=\Psi(\Hd,q+H^{\textsc{k,m}},L),$$
with $(\Hd,H^{\textsc{k,m}})$ a bivariate killed inhomogeneous subordinator, starting at $(x,0)$ under $\P_x$, with drift $(\frac{b^2}{2},0)$ and jump measure $\mud$, defined for all $a\in(0,\tau)$, $u\in(0,\tau-a)\times\{+\infty\}$ and $q\in\{0,1\}$ by :
$$\mud(a,\d u,\d q):=\frac{1}{W(a)}\delta_{(+\infty,0)}(\d u,\d q)+ \int_{(0,a)} \d x\,\L(x+\d u)\,\B{f(x+u)}(\d q)\,\frac{W(a-x)W(\tau-a-u)}{W(a)W(\tau)},$$
and $L:=\inf\{t\geq0,\ \Hd(t)=+\infty\}$ the killing time of $\Hd$.

\par\bigskip
We close this section by giving some explicit calculations in the cases where the limiting process $Z$ is either the standard Brownian motion, or an $\alpha$-stable Lévy process ($\alpha\in(1,2)$).
\paragraph*{Example 1 : The Brownian case}\quad

Consider the case where the population of $\Tree_n$ have exponential life spans with mean $1$. Then an appropriate rescaling of the JCCP of $\Tree_n$ leads in the limit to the standard Brownian motion. \\

We set :
$$\L_n(\d r)=e^{-r}\mathds1_{r\geq0}\d r\ \ \text{and}\ \ d_n=\frac{n^2}{2}.$$
Then, Assumption \HypLevy\ is satisfied : for all $\lambda\geq0$, we have $\tilde\psi_n(\lambda)=\frac{n}{\lambda+n}\frac{\lambda^2}{2}$, which converges to $\psi(\lambda)=\frac{\lambda^2}{2}$ as $n\to\infty$, and this implies the convergence in $\D(\R)$ of $\tilde Z_n$ towards the standard Brownian motion (see \cite[Th. VII.3.4]{JS}). Moreover, if we assume $\theta_n=\frac\beta n$ for some $\beta\in[0,1]$, Assumption \HypMutConst\ holds with $\theta=\frac\beta 2$.\\

The genealogical structure of this process (without mutations) and its asymptotic behaviour are studied by L. Popovic in \cite{Popovic}, and in particular, results taking into account a $\beta$-sampling of extinct individuals (each individual in the genealogy is recorded with a probability $\beta$) are provided. The following results are presented as a consequence of Theorem \ref{th_cv_PPP_Const} but can also be derived from \cite{Popovic}, since $\beta$-sampling can be directly interpreted as recording $1$-type birth events in the genealogy. \\

The distribution of $\Sigma$ is completely explicit. We know that $W(x)=2x$, and $\Ha(t)=\frac t2$ a.s. for all $t\geq0$. Note that the image by $\Ha$ of a Poisson process with parameter $\theta$ is itself a Poisson process, with parameter $2\theta$. As a consequence, if we denote by $((a_0,0),(a_1,1)...,(a_j,1))$ the ranked sequence of the atoms of the measure $\sigma$ appearing in Theorem \ref{th_cv_PPP_Const}, under $N''(\,\cdot\,\cap\,\sup\epsilon\in(0,\tau))$, conditional on $a_0$, $(a_1,...,a_j)$ is distributed as the sequence of jump times of a Poisson process with parameter $\beta$, restricted to $(0,a_0)$. \\

Besides, from the criticality of Brownian motion,  we have $N''(\,\cdot\,\times\{0,1\})=N'(\epsilon|_{[\chi,\zeta)}\in\,\cdot\,)$, and since an excursion of Brownian motion away from $0$ is such that $\chi=0$ or $\chi=\zeta$,
$$N''(\sigma\in\,\cdot\,,\ \sup\epsilon<\tau)=N'(\sigma\in\,\cdot\,,\ \sup\epsilon\in(0,\tau)).$$
Finally, we have 
$$N'(\Ha(L(T^0-))\in\d h,\  \sup\epsilon\in(0,\tau))=N'(\sup\epsilon\in\d h,\ \sup\epsilon\in(0,\tau))=\frac{\d h}{2h^2}\mathds1_{0<h<\tau}.$$
The measure $\Pi_1$ can then be expressed as follows :

$$\Pi_1=\int_0^\tau \frac{\d h}{2h^2} \int_{\mathcal M} \pi_{\beta,h}(\d \Theta) \mathds1_{\{\delta_{(h,0)}+\sum_{i\in I}\delta_{(\Theta_i,1)}\in\,\cdot\,\}},$$
where $\mathcal M$ denotes the space of point measures on $\R_+$, $\pi_{\beta,h}$ is the law of a Poisson process with parameter $\beta$ restricted to the interval $(0,h)$, and for any $\Theta\in\mathcal M$, $(\Theta_i)_{i\in I}$ denotes the sequence of jump times of $\Theta$.\\

In other words, in the limit the mutations on a lineage are distributed as an independent Poisson process with parameter $\beta$, stopped at an independent random time distributed as the depth of an excursion away from $0$, with depth lower than $\tau$. Note furthermore that simple calculations lead to $\Pi_1(B_{\geq1})=\infty$ and $\Pi(B_{\geq2})<\infty$ (using the notation introduced in Remark \ref{remark_mesure_infinie}) : the number of lineages carrying at least one mutation (resp. two mutations) is a.s. infinite (resp. finite).\\

Moreover, contrary to what is announced in Remark \ref{remark_first_lineage}, the loss of memory of the exponential distribution ensures here that there is no need to add extra assumptions to extend the result to the first lineage. In this case, the limiting object is then obtained by adding to $\Sigma$ a Dirac mass on $(0,\sigma_\tau+\delta_{(\tau,0)})$ where $\sigma_\tau$ is an independent Poisson process on $[0,\tau)$ with parameter $\beta$. \\

Finally, according to Remark \ref{remark_Brownien_B2}, these results are still valid for any choice of a sequence of functions $(f_n)$ and of a real number $\kappa$ satisfying \HypMut, replacing $\beta$ by $\kappa/2$.

\paragraph*{Example 2 : The stable case}\quad

Fix $\alpha\in(1,2)$ and set :
$$\L_n(\d r)=-\frac{r^{-\alpha-1}}{\Gamma(-\alpha)} \mathds1_{r>1}\d r\ \ \text{and}\ \ d_n=n^\alpha,$$

then we have for all $\lambda\geq0$, $\tilde\psi_n(\lambda)\to\lambda^\alpha$ which is the Laplace exponent of an $\alpha$-stable spectrally positive Lévy process and Assumption \HypLevy\ is satisfied. If we now set $\theta_n:=\beta/n^\alpha$ for some $\beta\in[0,1]$, Assumption \HypMutConst\ holds with $\theta=\beta$. \\

In this case we are able to characterize explicitly the inhomogeneous killed subordinator $\Hd$ defined in Proposition \ref{prop_mesure_de_sauts}. Indeed, we know that $Z$ has no Gaussian component, $\Lambda(\d z)=-\frac{z^{-\alpha-1}}{\Gamma(-\alpha)}\d z$, and $W(x)=\frac{x^{\alpha-1}}{\Gamma(\alpha)}$. Hence $\Hd$ has no drift and for all $a\in(0,\tau)$, $u\in(0,\tau-a)\times\{+\infty\}$, a simple calculation leads to 
$$\mud(a,\d u)=-\frac{u^{-\alpha-1}}{\Gamma(-\alpha)}\frac{au}{u+a} \left(\frac{\tau-a-u}{\tau}\right)^{\alpha-1}\d u+\frac{a^{\alpha-1}}{\Gamma(\alpha)}\delta_{+\infty}(\d u).$$

\par\bigskip

\section{Proofs of statements} \label{sec_Proof}
Proving our theorems first requires to give some preliminary results (Section \ref{sec_PR}), and in particular, the introduction of the marked ladder height process of $(\tZ_n,\tZ_n^\m)$ we described in the Introduction. The definition of this process and the convergence results we obtained in \cite{LPWM1} are reviewed in Section \ref{sec_cvMLHP}. Then Section \ref{sec_Proof_results} is devoted to the proof of the results stated in Section \ref{sec_Main_results}, relegating to Section \ref{sec_cv_MC} the proof of some technical result of convergence.

\subsection{Preliminary results}\label{sec_PR}
\subsubsection{Consequences of Assumption \HypLevy}
We state here some direct consequences of the convergence of $\tZ_n$ towards $Z$, which we prove in the appendix. Denote by $T_n^A$ the first entrance time of $\tZ_n$ in the Borel set $A$, and write $T_n^x$ for $T_n^{\{x\}}$. Recall that similar notation has been introduced in Section \ref{sec_Main_results} for the limiting process $Z$. Then Assumption \HypLevy\ leads to :

\begin{proposition}\label{prop_cv}
\begin{enumerate}[\upshape(i)]
 \item For all $x,y>0$, under $\P_0$, $\Txn$ (resp. $T_n^{(y,\infty)}$) converges in distribution to $\Tx$ (resp. $T^{(y,\infty)}$) as $n\to\infty$. \label{prop_cv_Tx}
 \item As $n\to\infty$, $\tilde\phi_n\to\phi$ uniformly on every compact set of $\R_+$, and in particular $\tilde\eta_n\to\eta$. \label{prop_cv_phi}
 \item As $n\to\infty$, $\tilde W_n\to W$ uniformly on $\R_+$, and $\tilde W'_n\to W'$ uniformly on every compact set of $\R_+^*$.  \label{prop_cv_W}
\end{enumerate}
\end{proposition}

\begin{remarque}
According to the remark after Lemma 8.2, and Exercise 8.4 in \cite{K}, in the infinite variation case the scale function of a Lévy process is differentiable on $\R_+^*$ with continuous derivative, and in the finite variation case, it has left and right derivatives on $\R_+^*$.
\end{remarque}

\subsubsection{Convergence of the marked ladder height process}\label{sec_cvMLHP}

In this section we define the marked ladder height process of $(\tZ_n,\tZ_n^\m)$, and recall the convergence results obtained for this process in \cite{LPWM1}.
\paragraph*{Local times at the supremum} \quad

We first need to specify local times at the supremum for the processes $\tZ_n$ and $Z$. We denote by $\F=(\F_t)_{t\geq0}$ (resp. $\F_n=(\F_{n,t})_{t\geq0}$) the natural filtration associated to $Z$ (resp. $\tZ_n$), that is for all $t\geq0$,
$$\F_t=\sigma\{Z_s,\ s\leq t\}\ (\text{resp. }\ \F_{n,t}=\sigma\{\tZ_n(s),\ s\leq t\})$$

For all $n\geq1$, let $(\tau_{n,i})_{i\geq0}$ be a sequence of i.i.d. random exponential variables, independent of $(\tZ_n)_{n\geq1}$, with parameter $\alpha_n:=\frac{d_n}{n}$. Then, according to \cite[Chapter IV]{B} (or see \cite[Section 2]{LPWM1}), we define for $\tZ_n$ a local time at the supremum as follows :
\[L_n(t):=\sum_{i=0}^{\ll_n(t)} \tau_{n,i},\]
where $\ll_n(t)$ represents the number of jumps of the supremum until time $t$. We denote by $\Linv_n$ the right-continuous inverse of $L_n$, and replace the filtration $\F_{n,t}$ with $\F_{n,t} \vee \sigma (L_n(s),\, s\leq t)$, so that $L_n$ (resp. $\Linv_n$) is adapted to $(\F_{n,t})$ (resp. to $(\F_{n,\Linv_n(t)})$).\\

We introduce the local time at the supremum $L$ for the infinite variation Lévy process $Z$ : it is defined up to a multiplicative constant, and we require that
\begin{equation} \label{formula_norm_local_time}
\E\bigg(\int_{(0,\infty)}e^{-t}\d L_t\bigg)=\phi(1),
\end{equation}
so that $L$ is uniquely determined. Finally, we denote by $\Linv$ its inverse.\\

\paragraph*{Excursions and mutations} \quad

From now on, we assume (unless otherwise specified) that $\tZ_n^\m(0)=0$ a.s. We denote by $(t,e_{n,t})_{t\geq0}$ the excursion process of $\tZ_n$ formed by the excursions from its past supremum, and $N_n$ its excursion measure, as defined in Section \ref{sec_SPLP}. 

We define for all $t\in[0,L_n(\infty))$
$$\xi_n:=\left\{ 
\begin{array}{l l}
 (t,e_{n,t}(\zeta),\Delta\Zbn(\Linv_n(t)))_{t \geq0} & \quad \text{if } \Linv_n(t-)<\Linv_n(t)\\
 \partial & \quad \text{else}\\
\end{array}
\right.,$$
where $\partial$ is an additional isolated point, and $e_{n,t}(\zeta)$ stands for $e_{n,t}(\zeta(e_{n,t}))$.

Here the fourth coordinate $\Delta \Zbn(\Linv_n(t))$ is $1$ or $0$ whether or not the jump of $\tZ_n$ at the right end point of the excursion interval indexed by $t$ is marked. Note that the set $\{\Linv_n(t)\}_{t\geq0}$ of these right end points is exactly the set of record times of $\tZ_n$.

\paragraph*{Marked ladder height process} \quad

Then according to \cite{LPWM1}, for $n\geq1$, we define the marked ladder height process $H_n=(\Han,\Hcn)$ of $(\tZ_n,\tZ_n^\m)$ as the (possibly killed) bivariate subordinator with no drift and whose jump point process is a.s. equal to $\xi_n$. Moreover, according to Proposition 3.2 in \cite{LPWM1}, $H_n$ has Lévy measure
\begin{equation} \label{mu_n}
\mu_n(\d y,\d q):=\int_0^\infty  \d x\; e^{-\teta_n x}\; \tL_n(x+\d y) \; \B{f_n(n(x+y))}(\d q),
\end{equation}
and is killed at rate $\kill_n=\frac{1}{\tilde W_n(\infty)}$ if $\tZ_n$ is subcritical.

Note that $\Han$ is in fact the ladder height process of $\tZ_n$, i.e. for all $t\geq0$, $\Han(t)=\bar \tZ_n(\Linv_n(t))$ a.s., where $\bar \tZ_n(t)$ denotes the current supremum of $\tZ_n$ at time $t$. The jumps of $\Hcn$ correspond, in the local time scale, to the marks occurring at record times of $\tZ_n$. Moreover, $\Hcn$ is a Poisson process with parameter $\lambda_n:=\mu_n(\R_+^*\times\un)$, so that the random time 
\begin{equation} \label{e_n}
\e_n:=\inf\{t\geq0,\ \Hcn(t)=1\}
\end{equation}
follows on $\{\e_n<L_n(\infty)\}$ an exponential distribution with parameter $\lambda_n$. 

\paragraph*{Convergence theorem for the marked ladder height process} 

We define 
$$\mu(\d u,\d q):=\int_0^\infty  \d x\; e^{-\eta x}\; \L(x+\d u) \; \B{f(x+u)}(\d q),$$
and
$$\mu^+(\d u):=\mu(\d u,\{0,1\})=\int_0^\infty  \d x\; e^{-\eta x}\; \L(x+\d u).$$

We recall here Theorem 4.1 of \cite{LPWM1} :
\begin{theoreme} \label{th_cv_H}
Under Assumption \HypMutConst, the sequence of bivariate subordinators $H_n=(\Han,\Hcn)$ converges weakly in law to a subordinator $H:=(\Ha,\Hc)$, which is killed at rate $\kill:=\frac1{W(\infty)}$ if $Z$ drifts to $-\infty$. Moreover, $\Ha$ and $\Hc$ are independent, $\Ha$ is a subordinator with drift $\frac{b^2}2$ and Lévy measure $\mu^+$, and $\Hc$ is a Poisson process with parameter $\theta$.
\\

 Under Assumption \HypMut, the sequence of bivariate subordinators $H_n=(\Han,\Hcn)$ converges weakly in law to a subordinator $H:=(\Ha,\Hc)$, which is killed at rate $\kill$ if $Z$ drifts to $-\infty$. Moreover, $H$ has drift $(\frac{b^2}2,0)$ and Lévy measure 
\upshape $\mu(\d u,\d q)+\rho\delta_0(\d u)\delta_1(\d q),$ \itshape where $\rho:=\kappa b^2$. 
\end{theoreme}

In particular, under Assumption \HypMut, if $Z$ has no Gaussian component, the limiting marked ladder height process is a pure jump bivariate subordinator with Lévy measure $\mu$. If $Z$ has a Gaussian component, the fact that the \og small jumps\fg\ of $\tZ_n$ generate the Gaussian part in the limit results in a drift for $\Ha$, and possibly additional independent marks that happen with constant rate in time, as under Assumption \HypMutConst. This rate is proportional to the Gaussian coefficient (provided that $\kappa\neq0$). Besides, note that as expected, $\Ha$ is distributed as the classical ladder height process of $Z$.
\par\bigskip
An easy adaptation of the proof of this theorem yields 
\begin{theoreme} \label{th_cv_H*}
 Let $H_n^*$ be a driftless subordinator on $\R_+$ with Lévy measure
 \upshape $$\mu_n^*(\d u):=\int_{(0,\infty)}  \d x\; e^{-\teta_n x}\; \tL_n(x+\d u) \; (1-f_n(n(x+u))).$$
  \itshape Then $H_n^*$ converges in distribution to a subordinator $H^*$ with drift $\frac{b^2}2$ and Lévy measure
\upshape $$\mu^*(\d u)=\int_{(0,\infty)}  \d x\; e^{-\eta x}\; \L(x+\d u) \; (1-f(x+u)).$$
\end{theoreme}

We denote by $\psi_n^*$ and $\psi^*$ the respective Laplace exponents of $H_n^*$ and $H^*$.

\begin{remarque}
 Under Assumption \HypMutConst, this theorem is not of interest, since $H_n^*$ (resp. $H^*$) is equal in law to $H_n^+$ (resp. $H^+$), and so the result is given by Theorem \ref{th_cv_H} with $f\equiv 0$.
\end{remarque}

Finally, we recall Theorem 5.1 of \cite{LPWM1} :
\begin{theoreme} \label{th_appendix}
 The following convergence in distribution holds in $\D(\R^4)$ as $n\to\infty$ :
$$(\tZ_n,L_n,\Han,\Hcn)\Rightarrow(Z,L,\Ha,\Hc).$$
\end{theoreme}

\subsection{Proof of main results} \label{sec_Proof_results}
The proof of Theorems \ref{th_cv_PPP_Const} and \ref{th_cv_PPP} is organized in four subsections. In the first one we describe the distribution of the point measures $\sigma\ni$ from a family of Markov chains. More precisely, we show that these point measures are i.i.d., and that for any $\eps\in(0,\tau)$, their restriction to $[\eps,\tau)\times\{0,1\}$ has the law of a point measure whose set of atoms forms a Markov chain $M\neps$, killed at some first entrance time. The second one deals with the construction of the limiting Markov chain $M_\eps$, and then with the proof of theorems themselves, in which we make use of the convergence in distribution of $(M\neps)_n$. The proof of the latter convergence is quite long and is gathered in the last two subsections.

\subsubsection{Distribution of the point measures $\sigma\ni$}

From the article \cite{ALContour} of A. Lambert, we know that there is a one-to-one correspondence between a splitting tree and its JCCP. In particular, properties linked to the lineage of the $i$-th extant individual at level $\tau$ are read from the $i$-th excursion under level $\tau$ of the truncated JCCP. Then using the invariance by time reversal of such excursions, and making use of the strong Markov property, we obtain the following proposition. Recall that we conditioned $\tTree_n$ on having $I_n$ extant individuals alive at $\tau$.

\begin{proposition}\label{prop_loi_sigma}
Fix $\eps\in(0,\tau)$, $n\geq1$, and let $\sigma\nieps$ denote the trace measure of $\sigma\ni$ on $[\eps,\tau)\times\{0,1\}$. Then we have :
\begin{enumerate}[(i)] 
\item The random measures $(\sigma\nieps)_{1\leq i< I_n}$ are i.i.d. 
\item There exists a Markov chain $M\neps$ with values in $[\eps,\tau)\times\{0,1\}$ such that with probability $1-p_{n,\eps}$, $\sigma\neps^{(1)}([\eps,\tau)\times\{0,1\})=0$, and with probability $p_{n,\eps}$, $\sigma\neps^{(1)}$ is distributed as 
$$\sum_{k=0}^{K\neps} \delta_{M\neps(k)},$$
where $p\neps:=\frac n{d_n}\frac{\frac1{\tilde W_n(\eps)}-\frac1{\tilde W_n(\tau)}}{1-\frac{\tilde W_n(0)}{\tilde W_n(\tau)}}$, and $K\neps:=\inf\{k\geq0,\ M\neps^2(k)=0\}$ ($M\neps^i$, $i\in\{1,2\}$, denoting the $i$-th coordinate of $M\neps$).
\end{enumerate}
\end{proposition}

The probability $p\neps$ has in fact to be interpreted as follows : we have
$$\P_0(T_n^{-\eps}<T_n^{(0,\infty)}<T_n^{-\tau})=\P_0(T_n^{-\eps}<T_n^{(0,\infty)})-\P_0(T_n^{-\tau}<T_n^{(0,\infty)})=\frac{\tilde W_n(0)}{\tilde W_n(\eps)}-\frac{\tilde W_n(0)}{\tilde W_n(\tau)},$$
and hence
$$p_n=\P_0(T_n^{-\eps}<T_n^{(0,\infty)} \sachant T_n^{(0,\infty)}<T_n^{-\tau}).$$

\paragraph{Construction of $M\neps$}\quad

We construct below the Markov chain $M\neps$ appearing in Proposition \ref{prop_loi_sigma}, and which will converge in distribution towards the Markov chain $M_\eps$ that appears in Theorem \ref{th_cv_PPP} (the proof of the latter point is the purpose of Sections \ref{sec_Proof_firstmutation} and \ref{sec_Proof_cvMeps}).\\

Recall that we defined in Section \ref{sec_cvMLHP} (formula \eqref{e_n}) the random variable $\e_n=\inf\{t\geq0,\ \Hcn(t)=1\}.$ We set for all $n\geq1$, $x>0$ and $u\geq0$ : 
\begin{align}\label{formula_def_num}
 \num_n(x,\d u) &:=\P_0\big(\Han(\e_n)\in \d u, \Linv_n(\e_n)<\Txn \sachant \Txn<T_n^{(\tau-x,\infty)}\big)\\
 \nud_n(x,\d u) &:=\P_0\big(\bar \tZ_n(\Txn)\in \d u, \Linv_n(\e_n)\geq\Txn \sachant \Txn<T_n^{(\tau-x,\infty)}\big),
\end{align}
 where the letters $\m$ and $\textsc{d}$ stand respectively for \og mutation\fg\ and \og death\fg. \\

We want to initialize the Markov chain $M\neps$ at the first $1$-type birth event that occurs below level $\tau-\eps$, when following the lineage backward in time. Then, conditional on $\tZ_n(0)=0$ and $T_n^{-\eps}<T_n^{(0,\infty)}<\infty$, we set
$$S_n:=\sup\{t\leq T_n^{(0,\infty)},\ \tZ_n(t)<-\eps \}\ \text{ and } (\Ups_n,\MUpsn):=(-\tZ_n(S_n-),\Delta \tZ_n^\m(S_n)).$$

\begin{figure}[!h]
\begin{center}
  \includegraphics[width=10cm]{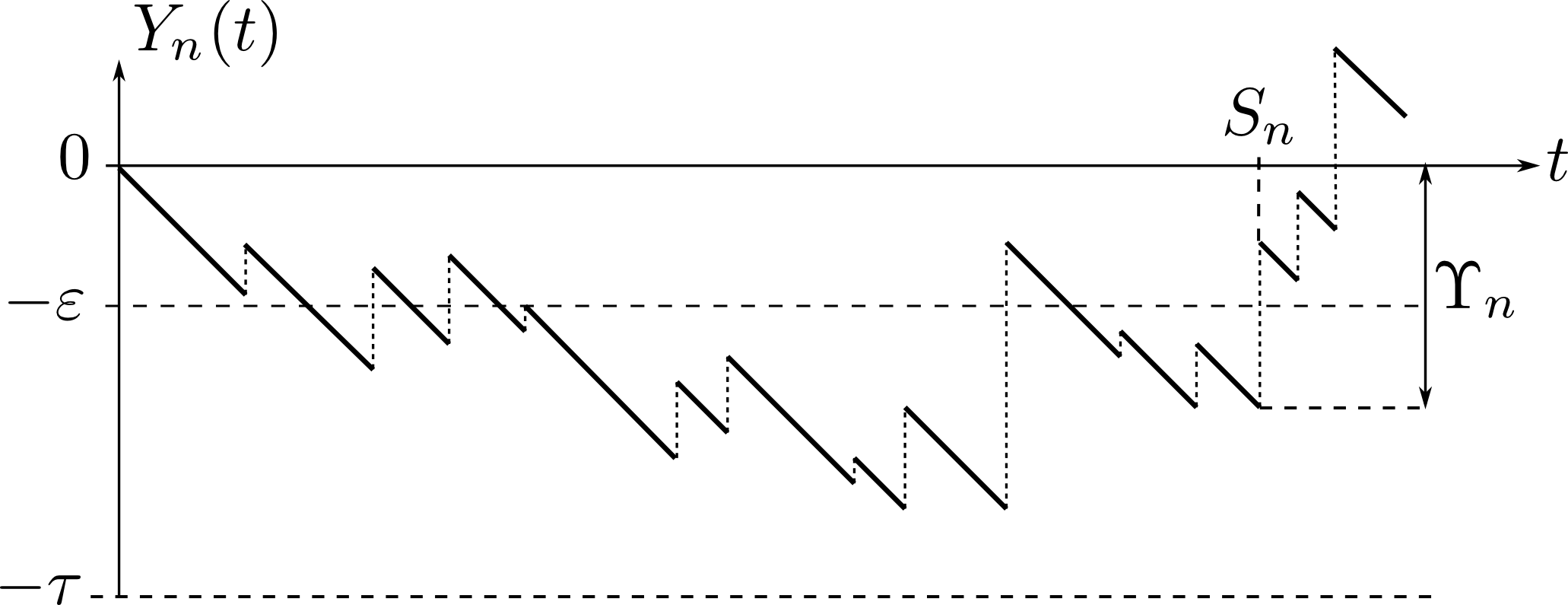}
 \caption{An excursion of $\tZ_n$ under level 0 and the random variables $S_n$ and $\Ups_n$.}
\label{figure_Ups}
\end{center}
 \end{figure}

Thereby if we consider an excursion of $\tZ_n$ away from $0$, $\Ups_n$ is the value of $\tZ_n$ before its last jump over level $-\eps$ (see Figure \ref{figure_Ups}), and $\MUpsn$ is the mark carried by this jump.

Finally we define, for all $(u,q)\in(\eps,\tau)\times\{0,1\}$,
$$\nuin(\d u,\d q):=\frac{1}{p'_{n,\eps}}\P_0(\Ups_n\in\d u,\ \MUpsn\in\d q,\ T_n^{-\eps}<T_n^{(0,\infty)}<T_n^{-\tau}),$$
where $p'_{n,\eps}:=\frac n{d_n}\left(\frac1{\tilde W_n(\eps)}-\frac1{\tilde W_n(\tau)}\right)=p\neps\left(1-\frac{\tilde W_n(0)}{\tilde W_n(\tau)}\right)$ is in fact equal to $\P_0(T_n^{-\eps}<T_n^{(0,\infty)}<T_n^{-\tau})$, and is therefore a normalizing constant such that $\nuin$ is a probability measure. \\

Then we consider the Markov chain $M\neps=(M\neps(k))_{k\in\Z_+}$ with values in $[\eps,\tau)\times\{0,1\}$, defined by :
\begin{itemize}
 \item For all $k\in\Z_+$, for all $u\geq0$, conditional on $M\neps(k)=(x,1)$,
\[\left\{
 \begin{array}{ll}
  M\neps(k+1)\in (x+\d u)\times\un & \text{ with probability } \ \num_n(x,\d u) \\
  M\neps(k+1)\in (x+\d u)\times\zero & \text{ with probability }\ \nud_n(x,\d u) .
 \end{array}\right.
\]
 \item For all $k\in\Z_+$, conditional on $M\neps(k)=(x,0)$, $M\neps(k+1)=(x,0)$ a.s.
 \item For all $u\in[\eps,\tau)$,
\[\left\{
 \begin{array}{l}
\P(M\neps(0)\in\d u\times\un)=\nuin(\d u\times\un)+\int_{[\eps,\tau)} \nuin(\d x\times\zero) \num_n(x,\d u-x)\\
\P(M\neps(0)\in\d u\times\zero)=\int_{[\eps,\tau)} \nuin(\d x,\zero)\nud_n(x,\d u-x)
\end{array}\right.\]
\end{itemize}

Recall that $K\neps=\inf\{k\geq0,\ M\neps^2(k)=0\}$. Then all the information we need is contained in $(M\neps(0),\ldots,M\neps(K\neps))$ : the $K\neps$ first values $M\neps(0)$ to $M\neps(K\neps-1)$, which have second coordinate $1$ a.s., describe the law of the successive levels where a mutation occurred on a lineage $i$ up to level $\tau-\eps$. The random variable $M\neps^1(K\neps)$ has the law of the coalescence time between the two consecutive extant individuals $i-1$ and $i$ at level $\tau$, and $M\neps^2(K\neps)=0$ a.s.

\paragraph{Proof of Proposition \ref{prop_loi_sigma}}\quad

We denote by $\Tree_{n,n\tau}$ the truncation of $\Tree_n$ up to level $n\tau$, and by $(Z_{n,n\tau},Z_{n,n\tau}^{\textsc m})$ the JCCP of $\Tree_{n,n\tau}$. We define 
$$(\tZntau(t),\tZntau^{\textsc m}(t))_{t\geq0}:=(\frac1n Z_{n,n\tau}(d_n t),Z_{n,n\tau}^{\textsc m}(d_n t))_{t\geq0},$$
 which is in fact, up to a rescaling of time, the JCCP of the rescaled marked splitting tree $\tTree_n$, truncated up to level $\tau$.

\begin{figure}[!h]
  \includegraphics[width=15.5cm]{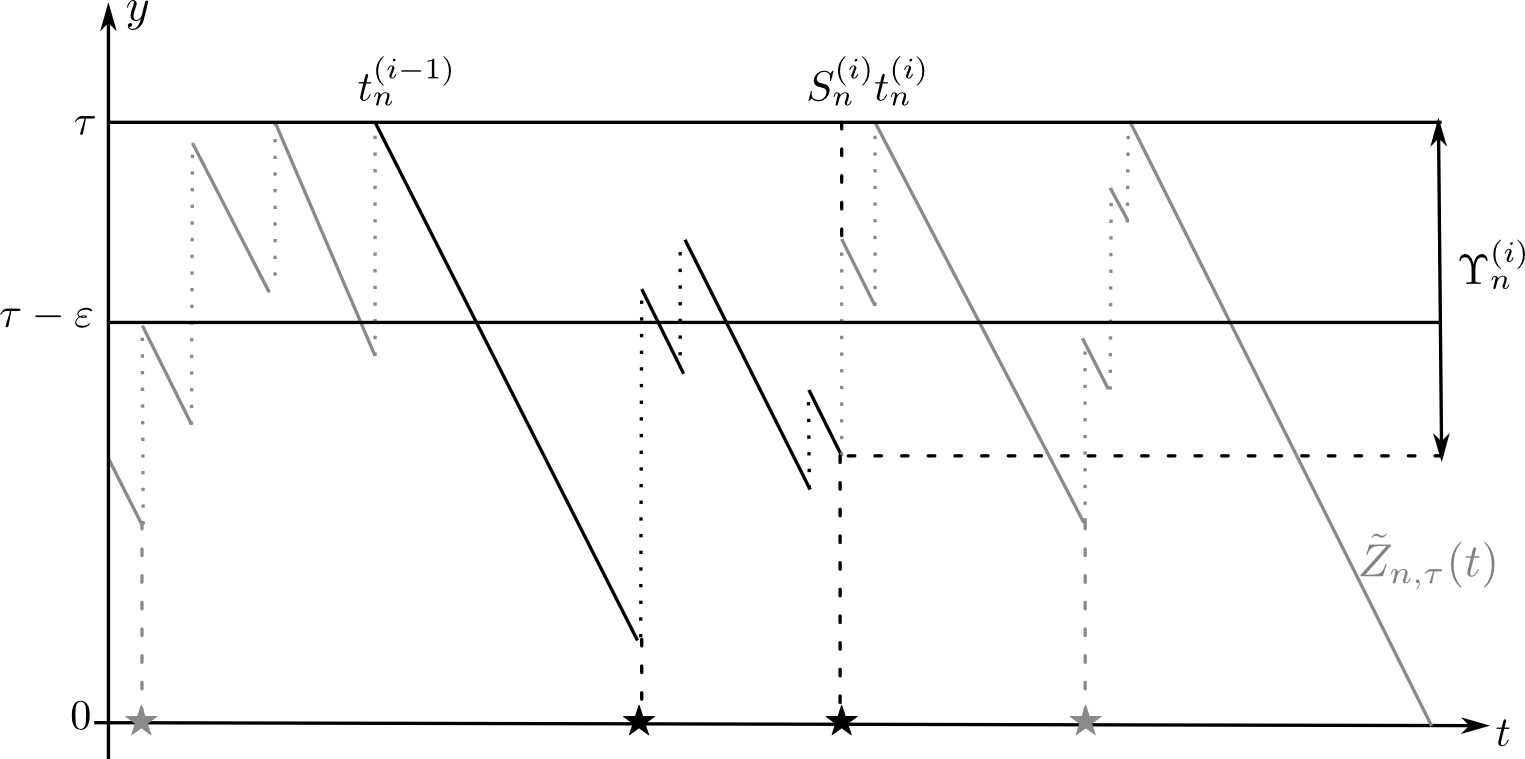}
 \caption{A representation of the (rescaled in time) JCCP $(\tZntau,\tZntau^\m)$ (where as before, $\tZntau^\m$ is represented by the sequence of its jump times, symbolized by stars on the horizontal axis). Here $\Ups_n^{(i)\textsc m}=1$. The reversed path $e\nieps$ can be read from the black path and black stars, reading the figure upside down and changing $y$ on the vertical axis into $\tau-y$.}
\label{figure_enieps}
 \end{figure}

The following lemma is a key tool for the analysis of the genealogy. See Figure \ref{figure_enieps} for graphical interpretation of some of the objects involved.

\begin{lemme}\label{lemme_exc_indep}
 Fix $n\geq1$ and $\eps>0$. Define : 
\begin{align*}
& t_n^{(0)}:=\inf\{t\geq 0,\ \tZntau(t)=\tau\},\ \text{ and for}\ i\in\N,\ t\ni:=\inf\{t>t_n^{(i-1)},\ \tZntau(t)=\tau\}, \\
& S\ni:=\sup\{t\in[t_n^{(i-1)},t\ni],\ \tZntau(t)<\tau-\eps \},\\
& \text{ and } (\Ups\ni,\Ups_n^{(i)\textsc m}):=(\tau-\tZntau(S\ni-),\Delta \tZntau^\m(S\ni)).
\end{align*}

Only the first $I_n$ values in the sequence $(t\ni)_{i\geq0}$ are finite, and the reversed killed paths 
$$\exc\ni:=\Big\{\big(\tau-\tZntau((t\ni-t)-),\tZntau^\m(t\ni)-\tZntau^\m((t\ni-t)-)\big),\ 0\leq t< t\ni-t_n^{(i-1)}\Big\},\ 1\leq i< I_n,$$ 
are i.i.d. Besides, defining for all $1\leq i< I_n$, 
$$\exc\nieps:=(\exc\ni(t),\ t\ni-S\ni \leq t< t\ni-t_n^{(i-1)}),$$
conditional on $(\Ups\ni,\Ups_n^{(i)\textsc m})=(x,q)$, $\exc\nieps$ has the law of $(\tZ_n,\tZ_n^\m)$, starting at $(x,q)$, conditioned on $\tZ_n$ hitting $0$ before $(\tau,\infty)$, and killed when $\tZ_n$ hits $0$.
\end{lemme}

\begin{demo}
 From Theorem 4.3 in \cite{ALContour} which characterizes the law of the JCCP of $\Tree_{n,\tau}$ (without marks), we deduce that the paths $\{\tZntau(t),\ t_n^{(i-1)}\leq t < t\ni\}$, $1\leq i<I_n$, are i.i.d and distributed as $\tZ_n$ starting from $\tau$, conditioned on hitting $(\tau,\infty)$ before $0$ and killed when hitting $(\tau,\infty)$. Adapting this property to our marked trees, the i.i.d. property of $\{\exc\ni,\ 1\leq i<I_n\}$ is now straightforward, and the second part of the lemma is then obtained either from an appeal to Proposition \ref{prop_Exc_mesure} along with the Markov property of $(\Han,\Hcn)$ at $L_n(T_n^{(\eps,\infty)}-)$, or using directly a time reversal argument at the last exit time $S\ni$ (see \cite[Th. 3.10]{Nagasawa}).
\end{demo}

\par\medskip

\begin{demopr}{\ref{prop_loi_sigma}}
To begin with, we deduce from \cite[Corollary 3.5]{ALContour} that for $1\leq i<I_n$, the set of levels at which birth events occurred on the $i$-th lineage is the set of the values taken by the future infimum of the rescaled JCCP between $t_n^{(i-1)}$ and $t\ni$, i.e. by the process $$j_n(t):=\underset{[t,t\ni]}\inf \tZntau,\ \ t_n^{(i-1)}\leq t\leq t\ni.$$ 
As a consequence, the subset of those levels corresponding to $1$-type birth events is a.s. equal to $\{j_n(t-),\ t\in J_n^\star \},$ where $J_n^\star$ is the set of jump times of $j_n$ (which are necessarily jump times of $\tZntau$) carrying a mark :
\[J_n^\star:=\{s\in(t_n^{(i-1)},t\ni],\ \Delta j_n(s)>0\ \text{and}\ \Delta \tZntau^\m(s)>0 \}.\]
Moreover from \cite[Theorem 3.4]{ALContour}, the coalescence time between lineage $i$ and lineage $i-1$ is given by $\tau-\inf_{[t_n^{(i-1)},t\ni]} \tZntau=\tau-j_n(t_n^{(i-1)})$. This yields $\sigma\ni=\delta_{(\tau-j_n(t_n^{(i-1)}),0)}+\sum_{t\in J_n^\star} \delta_{(\tau-j_n(t-),1)}$ a.s.\\

We are interested in the trace on $[\eps,\tau)\times\{0,1\}$ of $\sigma\ni$. From the preceding observations and using Lemma \ref{lemme_exc_indep}, we deduce the following : first, the point measures $\sigma\nieps$ are i.i.d. Second, since  $p\neps=\P_0(T_n^{-\eps}<T_n^{(0,\infty)} \sachant T_n^{(0,\infty)}<T_n^{-\tau})$, then with probability $1-p\neps$, the infimum of the excursion $\exc\ni$ is greater than $-\eps$, implying $\sigma\nieps([\eps,\tau)\times\{0,1\})=0$. Else with probability $p\neps$, the point measure $\sigma\nieps$ has at least one atom. \\

Conditional on $\sigma\nieps$ having at least one atom, we choose to order these atoms as in the definition of the space $\Mponct$, i.e. increasingly w.r.t. the first coordinate and decreasingly w.r.t. the second one. 
First note that the reversed future infimum $\big(\tau-j_n((t\ni-t)-),\ 0\leq t<t\ni-t_n^{(i-1)}\big)$ is a.s. equal to the running supremum of $\exc\ni$. Then, from Lemma \ref{lemme_exc_indep} and the first part of this proof, we deduce the following :
\begin{itemize}
 \item [$\circ$] Denote by $(a_0,q_0)$ the first atom of $\sigma\nieps$. Conditional on $(\Ups\ni,\Ups_n^{(i)\m})=(u,q)$, if $q=1$ we have $(a_0,q_0)=(u,1)$ a.s. If $q=0$, then $(a_0,q_0)\in u+\d v\times\un$ with probability $\num_n(u,\d v)$, and $(a_0,q_0)\in u+\d v\times\zero$ with probability $\nud_n(u,\d v)$. Consequently, $(a_0,q_0)$ is distributed as $M_\eps(0)$.
 \item [$\circ$] Now conditional on $(a_0,q_0)$, if $q_0=0$, then $\sigma\nieps$ has one unique atom. Now $M\neps^2(0)=0$ implies $K\neps=0$ a.s., so that we have as announced $\sigma\nieps\stackrel{(d)}{=}\sum_{k=0}^{K\neps} \delta_{M\neps(k)}$. Else if $q_0=1$, applying the strong Markov property to $(\Ha_n,\Hc_n)$ at $\e_n$, the next atom of $\sigma\nieps$ has the law of $M\neps(1)$ conditional on $M\neps(0)=(x,1)$. 
\end{itemize}

Finally, through a recursive application of the Markov property, stopped the first time an atom has second coordinate $0$, we obtain the announced equality in law.
\end{demopr}

\subsubsection{Limiting Markov chain}

Similarly as in the last subsection, for fixed $\eps\in(0,\tau)$ we define a Markov chain $M_\eps$, towards which the sequence $(M\neps)$ will converge in distribution. First define thanks to Theorem \ref{th_cv_H}
$$\e:=\inf\{t\geq0,\ \Hc(t)=1\}.$$
Note that as $\e_n$, $\e$ follows an exponential distribution, whose parameter $\lambda$ is equal to $\theta$ in the case of Assumption \HypMutConst, and to $\mu(\R_+^*,\un)+\rho$ in the case of Assumption \HypMut. \\

Then for all $x>0$, $u>0$ and $q\in\{0,1\}$, we set 
\begin{align*}
 \num(x,\d u) &:=\P_0(\Ha(\e)\in \d u, \Linv(\e)<\Tx \sachant \Tx<T^{(\tau-x,\infty)})\\
 \nud(x,\d u) &:=\P_0(\bar Z(\Tx)\in \d u, \Linv(\e)\geq\Tx \sachant \Tx<T^{(\tau-x,\infty)}).
\end{align*}

We now want to define $\nui$, the counterpart in the limit of the measure $\nuin$ defined at rank $n$. The limiting process $Z$ having infinite variation, this measure will necessarily be described in terms of excursions. 

Let $\epsilon\in\mathscr E'$ satisfying $-\inf\epsilon\in(\eps,\tau)$. We define
$$S^\eps:=\sup\{t\leq\zeta,\ \ \epsilon(t)<-\eps \}$$
the last exit time of $\epsilon$ away from $(-\infty,-\eps)$. We then set
$$\Ups^\eps(\epsilon):=-\epsilon(S^\eps-),\ \text{ and }\ \Delta\Ups^\eps(\epsilon):=\epsilon(S^\eps)-\epsilon(S^\eps-).$$

Recall that the bivariate Lévy process $(\tilde Z_n,\tilde Z_n^\m)$ does not converge in general, as observed in Remark \ref{remark_pas_de_cv_Zm}. Then defining a process of marked excursions in the limit is not possible, and for this reason we do not directly define the counterpart of the r.v. $\MUps_n$.\\

In the sequel, when $\eps$ is fixed, the notation $\Ups$ (resp. $\Delta\Ups$) stands for $\Ups^\eps(\epsilon)$ (resp. $\Delta\Ups^\eps(\epsilon)$). Then if we consider an excursion of $Z$ away from $0$ conditioned on hitting level $-\eps$, $\Ups$ is the value of $Z$ before its last entry into $(-\eps,\infty)$ (see Figure \ref{figure_Ups} for a representation in finite variation).

Finally, we define for all $(u,q)\in[\eps,\tau)\times\{0,1\}$ :
$$\nui(\d u,\d q):=\frac{1}{p_\eps}\int_{(u-\eps,\infty)}N'(\Ups\in\d u,\ \Delta\Ups\in\d v,\ -\inf\epsilon\in[\eps,\tau))\B{f(v)}(\d q),$$
where $p_\eps:=\frac{1}{W(\eps)}-\frac{1}{W(\tau)}$. According to lemma 9 in \cite{Obloj}, we have $p_\eps=N'(-\inf\epsilon\in(\eps,\tau))$, so that $\nui$ is a probability measure. \\

\par\bigskip

Next let $M_\eps=(M_\eps(k))_{k\in\Z_+}$ be the Markov chain with values in $[\eps,\tau)\times\{0,1\}$, defined by :
\begin{itemize}
 \item For all $k\in\Z_+$, for all $u\geq0$, conditional on $M_\eps(k)=(x,1)$,
\[\left\{
 \begin{array}{ll}
  M_\eps(k+1)\in (x+\d u)\times\un & \text{ with probability } \ \num(x,\d u) \\
  M_\eps(k+1)\in (x+\d u)\times\zero & \text{ with probability }\ \nud(x,\d u) 
 \end{array}\right.
\]
 \item For all $k\in\Z_+$, conditional on $M_\eps(k)=(x,0)$, $M_\eps(k+1)=(x,0)$ a.s.
 \item For all $u\in[\eps,\tau)$,
\[\left\{
 \begin{array}{l}
\P(M_\eps(0)\in\d u\times\un)=\nui(\d u\times\un)+\int_{[\eps,\tau)} \nui(\d x\times\zero) \num(x,\d u-x)\\
\P(M_\eps(0)\in\d u\times\zero)=\int_{[\eps,\tau)} \nui(\d x,\zero)\nud(x,\d u-x)
\end{array}\right.\]
\end{itemize}
The values $0$ and $1$ stand as earlier for the absence or presence of a mutation. \\
Let $K_\eps$ be defined as follows :
$$K_\eps:=\inf\{k\geq0,\ M_\eps^2(k)=0\}.$$

Under $\P_x(\point\sachant T^0<T^{(\tau,\infty)})$, the interval $[0,L(T^0))$ is a.s. finite, and $K_\eps+1$ is a.s. equal to the number of jumps of the counting process $\Hc$ on this interval, so that $K_\eps$ is a.s. finite.\\

The main argument needed for the proof of Theorems \ref{th_cv_PPP_Const} and \ref{th_cv_PPP} is given by the following proposition :
\begin{proposition} \label{prop_cv_M}
For all $k\geq0$, as $n\to\infty$, the $(k+1)$-tuple $(M\neps(0),...M\neps(k))$ converges in distribution towards $(M_\eps(0),...M_\eps(k))$.
\end{proposition}

 For now we admit this proposition and relegate its proof to Section \ref{sec_cv_MC}. We now have all the necessary ingredients to prove our main theorem.

\subsubsection{Proof of Theorems \ref{th_cv_PPP_Const} and \ref{th_cv_PPP}}

In this Section we assume that one of the two Assumptions \HypMutConst\ or \HypMut\ is satisfied. We first establish the convergence of $\Sigma_n$ towards a Poisson point measure with intensity $\Leb\otimes\Pi$, making use of the law of rare events for null arrays (see e.g. Theorem 16.18 in \cite{Kall}). The proof of Theorem \ref{th_cv_PPP}, which is valid both under \HypMutConst\ and \HypMut, will then consist in identifying the intensity measures $\Pi_2$ with the measure $\Pi$.\\

Our main objects of interest in this section are then the point measures $\Sigma_n=\sum_{i=1}^{I_n} \delta_{(\frac{in}{d_n},\sigma\ni)}$, where we recall that the random variables $\sigma\ni$ have values in the space $\Mponct$ defined in Section \ref{sec_topo}. Note that a measure $\delta_{(a_0,0)}+\sum_{i=1}^j\delta_{(a_i,1)}$ in $\Mponct$ is characterized by the set of first coordinates of its atoms $\{a_0,\ldots,a_j\}$. Then, if we denote by $B_{m,\eps}$ the subset of $\Mponct$ defined by
$$B_{m,\eps}=\{\sigma\in\Mponct,\ \sigma([\eps,\tau)\times\{0,1\})=m+1\},$$
the class $\mathscr C:=\{B_{m,\eps},\ m\in\Z_+,\ \eps\in(0,\tau)\}$ is a generating class for the trace $\sigma$-field on $\Mponct$.

\begin{proposition}\label{prop_formerth_cv_PPP}
The sequence $(\Sigma_n)$ converges in distribution towards a Poisson point measure $\Sigma$ on $[0,1]\times\Mponct$ with intensity measure 
\upshape $\Leb\otimes\Pi,$ \itshape
where $\Pi$ is a measure on $\Mponct$ characterized as follows : for all $m\in\Z_+$ and $\eps\in(0,\tau)$,
\[\Pi(B_{m,\eps})=p_\eps\P(K_\eps=m).
\]
\end{proposition}

\begin{demopr}{\ref{prop_formerth_cv_PPP}}
To begin with, we prove that as $n\to\infty$, $\E(\Sigma_n (B\times C))\to\E(\Sigma (B\times C))$ for any Borel set $B$ in $[0,1]$ and any measurable set $C$ of $\Mponct$. From Lemma \ref{lemme_exc_indep}, we know that the point measures $\sigma\ni$, $1\leq i< I_n$, are independent, yielding 
\begin{align*}
 \E\left(\Sigma_n (B\times C)\right)
 &= \sum_{i=1}^{I_n-1} \P\left(\frac{in}{d_n}\in B,\ \sigma\ni\in C\right) \\
 &= \left(\frac{d_n}{n}\P(\sigma_n^{(1)}\in C)\right)\left(\frac{n}{d_n}\sum_{i=1}^{I_n-1} \mathds1_{\frac{in}{d_n}\in B}\right)
\end{align*}

Recall that we assumed that $I_n\underset{n\to\infty}\sim \frac{d_n}n$. The second term in the right-hand side clearly converges in distribution towards Leb$(B)$, and it remains to prove the convergence of the first term. Now using a monotone class argument, it suffices to prove this convergence for sets $C$ in the class $\mathscr C$ defined above.\\

For all $\eps\in(0,\tau)$ and $m\in\Z_+$, we have by definition of $\sigma\neps^{(1)}$ and according to Proposition \ref{prop_loi_sigma} : 
$$\frac{d_n}n\P(\sigma_n^{(1)}\in B_{m,\eps})=\frac{d_n}n\P(\sigma\neps^{(1)}\in B_{m,\eps})= \frac{d_n}{n}p_{n,\eps}\P(K\neps=m) .$$

First for $m=0$, we then have 
\begin{align*}
\frac{d_n}{n}\P(\sigma_n^{(1)}\in B_{0,\eps})
&=\frac{d_n}{n}p_{n,\eps}\P(M\neps^2(0)=0)\\
&\underset{n\to\infty}\longrightarrow\ p_\eps \P(M_\eps^2(0)=0)\ =\ p_\eps\P(K_\eps=0),
\end{align*}
and for $m\geq1$,
\begin{align*}
 \frac{d_n}{n}\P(\sigma_n^{(1)}\in B_{m,\eps})
&=\frac{d_n}{n}p_{n,\eps}\P(M\neps^2(m-1)=1,\ M\neps^2(m)=0)\\
&\underset{n\to\infty}\longrightarrow\ p_\eps \P(M_\eps^2(m-1)=1,\ M_\eps^2(m)=0)\ =\ p_\eps\P(K_\eps=m),
\end{align*}
where the convergences are obtained from an appeal to Proposition \ref{prop_cv_M} and using the fact that $ \frac{d_n}{n} p_{n,\eps}\to p_\eps$.\\

Finally, we get for all $B,C\in\mathcal B([0,1])\times\mathscr C$ :
$$\E\left(\Sigma_n (B\times C)\right)\underset{n\to\infty}\to \Pi(B \times C).$$

The point measures $(\Sigma_n)$ form a null array of simple point measures on $[0,1]\times\Mponct$, therefore, from the conclusion above, the theorem is a straightforward consequence of Theorem 16.18 in \cite{Kall}.
\end{demopr}

\par\bigskip
The following lemma is the last step preceding the proof of Theorems \ref{th_cv_PPP_Const} and \ref{th_cv_PPP}. 
For $i\geq1$, we define the sequence $(\e_i)_{i\geq0}$ as follows : first set $\e_0=0$. Then, for $i\geq1$, $\e_i$ denotes the $i$-th jump time of $\Hc$ if it exists, or is else set to $+\infty$. Note that $\e_1$ is in fact equal to $\e$ a.s. We then define $J:=\sup\{i\geq0,\ \e_i< L(T^0)\}$, which is in particular finite a.s. on $L(T^0)<\infty$.
\begin{lemme}\label{lemme_J}
For all $m\in\Z_+$ and $\eps\in(0,\tau)$ we have \upshape
\begin{equation} \label{eq3}
  \P(K_\eps=m)=\int_{[\eps,\tau)\times\{0,1\}} \nui(\d x,\d q)\P_x(J=m-q\sachant T^0<T^{(\tau,\infty)}).
\end{equation}
\end{lemme}

\begin{remarque}\label{remark_J}
Let $\sigma$ be defined as in Theorem \ref{th_cv_PPP}. Then for $x,\eps\in(0,\tau)$, $m\in\Z_+$, if $x\geq\eps$ then $\P_x(J=m \sachant T^0<T^{(\tau,\infty)})$ is in fact equal to $\P_x(\sigma\in B_{m,\eps}\sachant T^0<T^{(\tau,\infty)})$.
\end{remarque}

\begin{demo}
Fix $\eps\in(0,\tau)$. First note that for all $x\in[\eps,\tau)$,
\begin{equation}\label{eqJ0}
\P_x(J=0\sachant T^0<T^{(\tau,\infty)})
=\P_x(\e_1\geq L(T^0)\sachant T^0<T^{(\tau,\infty)})
=\nud(x,[\eps,\tau)),
\end{equation}
and
\begin{align*}
 \P_x&(J=1\sachant T^0<T^{(\tau,\infty)})\\
=&\,\P_x(\e_1<L(T^0),\ \e_2\geq L(T^0)\sachant T^0<T^{(\tau,\infty)})\\
=&\int_{[0,\tau-x)}\P_x(\e_1<L(T^0),\ \e_2\geq L(T^0),\ \Ha(\e_1)\in x+\d u,\ T^0<T^{(\tau,\infty)})/\P_x(T^0<T^{(\tau,\infty)})\\
=&\int_{0,\tau-x)}\P_x(\e_1<L(T^0),\ \Ha(\e_1)\in x+\d u)\P_{x+u}(\e_1\geq L(T^0),\ T^0<T^{(\tau,\infty)})/\P_x(T^0<T^{(\tau,\infty)})\\
=&\int_{[0,\tau-x)}\frac{\P_x(\e_1<L(T^0),\ \Ha(\e_1)\in x+\d u,\ T^0<T^{(\tau,\infty)})}{\P_x(T^0<T^{(\tau,\infty)})}\ \frac{\P_{x+u}(\e_1\geq L(T^0),\ T^0<T^{(\tau,\infty)})}{\P_{x+u}(T^0<T^{(\tau,\infty)})}\\
=&\int_{[0,\tau-x)}\P_x(\e_1<L(T^0),\ \Ha(\e_1)\in x+\d u \sachant T^0<T^{(\tau,\infty)})\P_{x+u}(\e_1\geq L(T^0)\sachant T^0<T^{(\tau,\infty)}),
\end{align*}
where in the third equality we applied the Markov property to $(\Ha,\Hc)$ at the stopping time $\e_1$. We omit for now to justify properly this application of the Markov property : details on filtrations and stopping times are provided in Section \ref{sec_Proof_firstmutation} (see Proposition \ref{prop_StoppingTimes}). 
Finally, this gives :
\begin{equation}
 \label{eqJ1} \P_x(J=1\sachant T^0<T^{(\tau,\infty)})=\int_{u\in[0,\tau-x)}\num(x,\d u)\nud(x+u,[\eps,\tau)).
\end{equation}

We first show \eqref{eq3} for $m=0$. Since $J\geq0$ a.s., from \eqref{eqJ0} we have
\begin{align*}
 \int_{[\eps,\tau)\times\{0,1\}} \nui(\d x,\d q)\P_x(J=-q\sachant T^0<T^{(\tau,\infty)})&= \int_{[\eps,\tau)} \nui(\d x,\{0\})\nud(x,[\eps,\tau))\\
&= \P(M^2_\eps(0)=0)= \P(K_\eps=0).
\end{align*}

Similarly we prove \eqref{eq3} for $m=1$, using \eqref{eqJ0} and \eqref{eqJ1} in the second equality : 
\begin{align*}
& \int_{[\eps,\tau)\times\{0,1\}} \nui(\d x,\d q)\P_x(J=1-q\sachant T^0<T^{(\tau,\infty)})\\
&= \int_{[\eps,\tau)}  \nui(\d x,\zero)\P_x(J=1\sachant T^0<T^{(\tau,\infty)})+\int_{[\eps,\tau)} \nui(\d v,\un)\P_v(J=0\sachant T^0<T^{(\tau,\infty)})\\
&= \int_{[\eps,\tau)}  \nui(\d x,\zero)\left(\int_{u\in[0,\tau-x)}\num(x,\d u)\nud(x+u,[\eps,\tau))\right)+\int_{[\eps,\tau)}  \nui(\d v,\un)\nud(v,[\eps,\tau))\\
&=\int_{v\in[\eps,\tau)}  \left( \nui(\d v,\un)+\int_{x\in[\eps,v)}\nui(\d x,\zero)\num(x,\d v-x)\right)\nud(v,[\eps,\tau))\\
&=\int_{[\eps,\tau)\times\{0,1\}}\P(M_\eps(0)\in\d u\times\{1\})\P(M^2_\eps(1)=0\sachant M_\eps(0)=(u,1))\\
&=\P(K_\eps=1).
\end{align*}
It is then clear by induction on $m$ that \eqref{eq3} is true for all $m\in\Z_+$, which ends the proof.
\end{demo}

\par\bigskip
In the proof below, we use Proposition \ref{prop_formerth_cv_PPP} and Lemma \ref{lemme_J} to deduce Theorem \ref{th_cv_PPP}, which is in fact also valid both under \HypMutConst\ and \HypMut. Theorem \ref{th_cv_PPP_Const} is then simply a consequence of  Theorem \ref{th_cv_PPP}, using the independence between $\Ha$ and $\Hc$ that arises under Assumption \HypMutConst. \\

\begin{demoth}{\ref{th_cv_PPP}}
Fix $m\in\Z_+$ and $\eps\in(0,\tau)$. First, from Proposition \ref{prop_formerth_cv_PPP}, along with Lemma \ref{lemme_J} and Remark \ref{remark_J}, we deduce
\begin{equation} \label{eqPiBmeps}
 \Pi(B_{m,\eps})=p_\eps \int_{[\eps,\tau)\times\{0,1\}} \nui(\d x,\d q)\P_x(\sigma\in B_{m-q,\eps}\sachant T^0<T^{(\tau,\infty)}).
\end{equation}

We now want to prove that $\Pi$ and $\Pi_2$ coincide on the generating class $\mathscr C$, using \eqref{eqPiBmeps}. We denote by $\sigma$ the point measure $\Psi(\Ha,\epsilon^\m+\Hc,L(T^0))$ that appears in the statement of the theorem, and we consider 
$$\Pi_2(B_{m,\eps})=N''(\sigma\in B_{m,\eps},\ \sup\epsilon<\tau).$$ 
Recall first that any point measure belonging to $B_{m,\eps}$ necessarily has at least one atom with first coordinate greater than $\eps$. Using the (slightly abusive) notation $T^{(\eps,\infty)}$ for the first entrance time in $(\eps,\infty)$ of an excursion $\epsilon\in\mathscr E''$, we apply the Markov property to $(\Ha,\Hc)$ at $L(T^{(\eps,\infty)})$ : recall that $\Ha(L(T^{(\eps,\infty)}))=Z(T^{(\eps,\infty)})$, and that $\sigma$ might have an atom coming from a jump of $\Hc$ at $L(T^{(\eps,\infty)})$. Conditional on $\Delta\epsilon(T^{(\eps,\infty)})= v$, this occurs with probability $f(v)$. Again, see Section \ref{sec_Proof_firstmutation} for details about filtrations and stopping times. This gives :
\begin{align*}
&N''(\sigma\in B_{m,\eps},\ \sup\epsilon<\tau)\\
&= \int_{[\eps,\tau)} \int_{[u-\eps,\infty)}N''(\sigma\in B_{m,\eps},\ \sup\epsilon\in[\eps,\tau),\ \epsilon(T^{(\eps,\infty)})\in\d u,\ \Delta\epsilon(T^{(\eps,\infty)})\in\d v)\\
&= \int_{[\eps,\tau)\times\{0,1\}} \int_{[u-\eps,\infty)} N''(\epsilon(T^{(\eps,\infty)})\in\d u,\ \Delta\epsilon(T^{(\eps,\infty)})\in\d v)\ \B{f(v)}(\d q)\\
&\qquad\qquad\qquad\qquad\qquad\qquad\qquad\qquad\qquad\qquad\qquad\qquad\qquad\times\P_u(\sigma\in B_{m-q,\eps},\ T^0<T^{(\tau,\infty)})\\
&= \int_{[\eps,\tau)\times\{0,1\}} \int_{[u-\eps,\infty)} N''(\epsilon(T^{(\eps,\infty)})\in\d u,\ \Delta\epsilon(T^{(\eps,\infty)})\in\d v)/\P_u(T^0<T^{(\tau,\infty)})\ \B{f(v)}(\d q)\\
&\qquad\qquad\qquad\qquad\qquad\qquad\qquad\qquad\qquad\times\P_u(\sigma\in B_{m-q,\eps},\ T^0<T^{(\tau,\infty)})\P_u(T^0<T^{(\tau,\infty)})\\
&= \int_{[\eps,\tau)\times\{0,1\}} \int_{[u-\eps,\infty)} N''(\epsilon(T^{(\eps,\infty)})\in\d u,\ \Delta\epsilon(T^{(\eps,\infty)})\in\d v,\ \sup\epsilon\in[\eps,\tau))\ \B{f(v)}(\d q)\\
&\qquad\qquad\qquad\qquad\qquad\qquad\qquad\qquad\qquad\qquad\qquad\qquad\qquad\times\P_u(\sigma\in B_{m-q,\eps}\sachant T^0<T^{(\tau,\infty)}).
\end{align*}

 Now from the definition of $N''$, we know that 
$$ N''(\epsilon(T^{(\eps,\infty)})\in\d u,\ \Delta\epsilon(T^{(\eps,\infty)})\in\d v,\ \sup\epsilon\in[\eps,\tau))= N'(\Upsilon\in\d u,\ \Delta\Upsilon \in\d v,\ -\inf\epsilon\in[\eps,\tau)),$$ which entails

$$N''(\sigma\in B_{m,\eps},\ \sup\epsilon<\tau)=p_\eps \int_{[\eps,\tau)\times\{0,1\}}  \nui(\d u,\d q)\P_u(\sigma\in B_{m-q,\eps}\sachant T^0<T^{(\tau,\infty)}).$$
This equality, along with  \eqref{eqPiBmeps}, leads to the expected result.
\end{demoth}

\par\bigskip
\begin{demoth}{\ref{th_cv_PPP_Const}}
As announced, the latter proof is also valid under \HypMutConst, in which case $\Hc$ is independent from $\Ha$, and is a Poisson process with parameter $\theta$. Moreover, $f\equiv0$ implies $N''(\epsilon^\m=1)=0$. Thus Theorem \ref{th_cv_PPP_Const} can be directly deduced from Theorem \ref{th_cv_PPP}.
\end{demoth}
\par\bigskip

\subsubsection{Proof of Proposition \ref{prop_mesure_de_sauts}}

Finally, we prove here Proposition \ref{prop_mesure_de_sauts}. The counterpart of this proposition under \HypMut\ (stated in the second paragraph of Section \ref{sec_Main_results}) can be established by an easy adaptation of the upcoming proof.

\par\bigskip
\begin{demopr}{\ref{prop_mesure_de_sauts}}
Fix $x\in(0,\tau)$. Consider the process $\Ha$ under $\P_x(\,\cdot\,\cap\, \{T^0<T^{(\tau,\infty)}\})$, killed at $L(T^0)$. This process is an inhomogeneous killed subordinator, with jump measure denoted by $\nuk$. Hereafter we prove that $\nuk$ and $\mud$ coincide.\\

 Let $F$ be a nonnegative continuous $\F_{\Linv}$-measurable function on $\R_+\times(\R_+\cup\{+\infty\})$, and $U$ a $\F_{\Linv}$-predictable process. Recalling that $L(T^0)$ is a $\F_{\Linv}$-stopping time, we have by compensation formula for any fixed $t>0$ :
\begin{multline}\label{eq5}
 \E_x \left(\sum_{0<r\leq t\wedge L(T^0)}  \left(\mathds1_{\Delta  \Ha_r>0}\, U_r\, F(\Ha_{r-},\Delta \Ha_r )\right) ,\ T^0<T^{(\tau,\infty)} \right)\\
=\ \E_x\left(\int_0^{t\wedge L(T^0)} \d s U_s \int_{(0,+\infty]} F(\Ha_s,z) \nuk(\Ha_s,\d z)\right)
\end{multline}
where $\Delta \Ha_r:=+\infty$ if $\Ha_r=+\infty$, and $\Delta \Ha_r:=\Ha_r-\Ha_{r-}$ otherwise. \\ 

On the other hand, we have:
\begin{align*}
 &\E_x\Bigg(\sum_{0<r\leq t\wedge L(T^0)}  \left(\mathds1_{\Delta  \Ha_r>0} \, U_r \,F(\Ha_{r-},\Delta \Ha_r ) \right),\ T^0<T^{(\tau,\infty)}\Bigg)\\
  &=\E_x\Bigg(\sum_{0<r\leq t} \E\Big( \,U_r\, F(\Ha_{r-},\Delta \Ha_r),\ \Delta  \Ha_r>0,\ r<L(T^0),\ T^0<T^{(\tau,\infty)} \sachant \F_{L^{-1}(r)}\Big)\\
 &\qquad\qquad\qquad\qquad\qquad\qquad\qquad\quad +\,U_r\, F(\Ha_{r-},\Delta \Ha_r),\ \Delta  \Ha_r>0,\ r=L(T^0),\ T^0<T^{(\tau,\infty)}  \Bigg)\\
  &=\E_x\Bigg(\sum_{0<r\leq t\wedge L(T^0)} \,U_r\, F(\Ha_{r-},\Delta \Ha_r) \,\left(\mathds1_{\Ha_{r-}<\Ha_r<\tau}\,\P_{\Ha_r}(T^0<T^{(\tau,\infty)})+\mathds1_{\Ha_{r-}<\tau,\ \Delta \Ha_r=+\infty}\right)\Bigg),
\end{align*}
using on the one hand the $\F_{\Linv(r)}$-measurability of every term but $\mathds1_{T^0<T^{(\tau,\infty)}}$ in the conditional expectation and the Markov property at time $r$, and on the other hand the fact that $\{r=L(T^0),\ T^0<T^{(\tau,\infty)}\}$ and $\{\Ha_{r-}<\tau,\ \Delta \Ha_r=+\infty\}$ coincide under $\E_x$.

We now express the sum in the right hand side in terms of excursions.
\begin{align*}
 &\E_x\left(\sum_{0<r\leq t\wedge L(T^0)}  \left(\mathds1_{\Delta  \Ha_r>0} \, U_r \,F(\Ha_{r-},\Delta \Ha_r ) \right),\ T^0<T^{(\tau,\infty)}\right)\\
   &=\E_x\Bigg(\sum_{0\leq g< \Linv(t)\wedge T^0} \,U_{L(g)}\, F(\Ha_{L(g)-},e_g(\zeta)) \mathds1_{\{-\inf e_g<\Ha_{L(g)-}<\tau-e_g(\zeta)\}}\, \P_{\Ha_{L(g)}}(T^0<T^{(\tau,\infty)})\\
   &\qquad\qquad\qquad\qquad\qquad\qquad\qquad\qquad\qquad\qquad\qquad\quad+U_{L(g)}\, F(\Ha_{L(g)-},+\infty)\,\mathds1_{\{-\inf e_g\geq\Ha_{L(g)}\}} \Bigg),
\end{align*}
where the sum in the right-hand side is taken over all the left-end points of excursions intervals.\\
Then by compensation formula,
\begin{multline}\label{eq6}
 \E_x\left(\sum_{0<r\leq t\wedge L(T^0)}  \left(\mathds1_{\Delta  \Ha_r>0} \, U_r \,F(\Ha_{r-},\Delta \Ha_r ) \right),\ T^0<T^{(\tau,\infty)}\right)\\
=\E_x\Bigg(\int_0^{t\wedge L(T^0)}\d s\, U_s\Bigg(\int_{(0,\tau-\Ha_s)} F(\Ha_s,z) \,\P_{\Ha_s+z}(T^0<T^{(\tau,\infty)})N(\epsilon(\zeta)\in\d z,\,-\inf\epsilon<\Ha_s)\\
+\,F(\Ha_s,+\infty)N(-\inf\epsilon\geq\Ha_s)\Bigg)\Bigg) 
\end{multline}
Finally from \eqref{eq5} and \eqref{eq6} we deduce that for all $a\in(0,\tau)$, $z\in(0,\infty]$,
$$\nuk(a,\d z)=\mathds1_{z<\tau-a}\P_{a+z}(T^0<T^{(\tau,\infty)}) N(\epsilon(\zeta)\in\d z,\,-\inf\epsilon<a) + N(-\inf\epsilon\geq a)\delta_{+\infty}(\d z),$$
which yields, using Proposition \ref{prop_Exc_mesure} and the fact that $N(-\inf\epsilon\geq a)=\frac{1}{W(a)}$,
$$\nuk(a,\d z)=\frac{W(\tau-a-z)}{W(\tau)} \int_0^a \d x \frac{W(a-x)}{W(a)}\L(x+\d z) + \frac{1}{W(a)}\delta_{+\infty}(\d z)=\mud(a,\d z).$$

From this result we deduce that under $\P_x$, $\Hd$ has the law of $\Ha$ under $\P_x(\,\cdot\,\cap \{T^0<T^{(\tau,\infty)}\})$, killed at $L(T^0)$, which finishes the proof.
\end{demopr}

\subsection{Convergence of the Markov chains}  \label{sec_cv_MC}
\subsubsection{Weak convergence of $\num_n$ towards $\num$ and characterization of these measures} \label{sec_Proof_firstmutation}

Before proving the convergence in law of $M\neps$ to $M_\eps$, we show in this subsection that the sequence of measure $(\num_n)$ converges weakly towards $\num$. 
Recall that $\num_n(x,\cdot)$ is the law of the amount of time elapsed between two mutations conditional on the latest one to have happened at level $\tau-x$ :

$$\num_n(x,\d u) :=\P_0\big(\Han(\e_n)\in \d u, \Linv_n(\e_n)<\Txn \sachant \Txn<T_n^{(\tau-x,\infty)}\big).$$

The announced weak convergence of $\num_n$ towards $\num$ is contained in the following result, which also gives an expression of these measures. Recall that in case $Z$ drifts to $-\infty$, we denoted by $\kill=\frac1{W(\infty)}$ the killing rate of $(\Ha,\Hc)$. If $Z$ does not drift to $-\infty$, we set $\kill=0$. 

\begin{theoreme} \label{th_cv_loi_des_sauts}
For all $z,y$ in $\R_+$ such that $z+y\leq\tau-x$, the measure
\upshape $$\P(\Han(\en-)\in\d z,\Delta\Han(\en)\in\d y,\Linv_n(\en)<\Txn \sachant \Txn<T_n^{(\tau-x,\infty)})$$
\itshape converges weakly towards
 \upshape  $$\P(\Ha(\e-)\in\d z,\Delta\Ha(\e)\in\d y,\Linv(\e)<\Tx \sachant \Tx<T^{(\tau-x,\infty)}) .$$
Besides, we have
\begin{multline} \label{formula_cv}
\P(\Ha(\e-)\in\d z,\ \Delta\Ha(\e)\in\d y,\ \Linv(\e)<\Tx<T^{(\tau-x,\infty)}) \\
=\Bigg\{\tilde\mu(\d y,\un)\Bigg[ U_*^{(\lambda+\kill)}(\d z)-\int_{[0,z)}\pi(\d a)\int_{[a,z)} U_*^{(\lambda+\kill)}(\d z-b) g^x(a,\zero,\d b-a)\Bigg]\\
\qquad\qquad\qquad\qquad\qquad\qquad\qquad\qquad-\pi(\d z)g^x(z,\un,\d y)
\Bigg\}\frac{W(\tau-x-z-y)}{W(x)} ,
\end{multline}
\itshape where
\begin{itemize}
 \item $\tilde\mu$ is the Lévy measure of $(\Ha,\Hc)$, yielding \upshape$\tilde\mu(\d y,\un)=\theta\delta_0(\d y)$\itshape\ under \HypMutConst, and \\ \upshape$\tilde\mu(\d y,\un)=\mu(\d y,\un)+\rho\delta_0(\d y)$\itshape\ under \HypMut.
 \item \itshape $U_*^{(l)}$ is the $l$-resolvent measure of the subordinator $H^*$ defined in Theorem \ref{th_cv_H*}, that is
\upshape $$U_*^{(l)}(\d z):=\int_{(0,\infty)} e^{-l t} \P(H^*(t)\in\d z) \d t,$$
 \item \itshape $\pi$ is a finite measure defined by
\upshape $$\pi(\d z):=\P(\Ha(L(\Tx)-)\in\d z,L(\Tx)\leq\e)$$
 \item \itshape and finally, \upshape
\begin{multline}
 \label{formula_g}
g^x(a,\d q,\d v)=\frac{b^2}2 (W'(x+a)-\eta W(x+a))\delta_0(\d v)\delta_0(\d q)\\
+\int_{(0,\infty)} (e^{-\eta u}W(x+a)-W(x+a-u))\;\B{f(u+v)}(\d q)\;\L(u+\d v)\d u.
\end{multline}
\end{itemize}
\end{theoreme}

Recall that $\lambda$ is the parameter of the exponential variable $\e$, and is equal to $\theta$ (resp. $\mu(\R_+^*,\{1\})+\rho$) under Assumption \HypMutConst\ (resp. \HypMut).

\begin{remarque}\label{remark_gx}
 In the case of Assumption \HypMutConst\ several simplifications can be made : we know that $\Ha$ and $\Hc$ are independent, and $f\equiv0$. Then  the processes $\Ha$ and $H^*$ are equal in law in $\D(\R_+)$, and further  \upshape $U_*^{(\lambda)}(\,\cdot\,)=\P(\Ha(\e)\in\,\cdot\,)$\itshape. Second, we have \upshape$g^x(a,\un,\d v)=0$ \itshape for all $x>0$ and $a,v\geq0$, and from \cite[(8.29)]{K}, we see that \upshape$g^x(a,\zero,\d v)=\P_{-(x+a)}(Z(T^{(0,\infty)})\in\d v)$\itshape. Finally, \eqref{formula_cv} yields
\upshape
\begin{multline} \label{formula_cv_Const}
\P(\Ha(\e-)\in\d z,\ \Delta\Ha(\e)\in\d y,\ \Linv(\e)<\Tx<T^{(\tau-x,\infty)}) \\
=\theta\delta_0(\d y)\Bigg[ U_*^{(\lambda)}(\d z)-\int_{[0,z)}\pi(\d a)\int_{[a,z)} U_*^{(\lambda)}(\d z-b) g^x(a,\zero,\d b-a)\Bigg]\frac{W(\tau-x-z-y)}{W(\tau)},
\end{multline}\itshape
which will be proven along with the theorem.
\end{remarque}

\begin{remarque}
 The measure $\pi$ is not explicit, and under Assumption \HypMut\ the random variable \upshape$\e$\itshape\ is not independent of $\Ha$ and $L(\Tx)$. However we can give another interpretation of $\pi$ in terms of a Poisson point measure : define similarly as in Section \ref{sec_cvMLHP}, for all $t\geq0$,
$$\xi(t):=\left\{ 
\begin{array}{l l}
 (e_t(\zeta),-\inf_{(0,\zeta)} e_t, \Delta \Hc(t))_{t \geq0} & \quad \text{if } \Linv(t-)<\Linv(t)\\
 \partial & \quad \text{else}\\
\end{array}
\right.,$$
where $\partial$ is an additional isolated point, and $(t,e_t)_{t\geq0}$ the excursion process of $Z$ (excursions from the past supremum). Then $(t,\xi(t))_{t\geq0}$ is a Poisson point process with values in $\R_+\times(\R_+^*)^2\times\{0,1\}$. Denote by $m$ its intensity measure, and by $\xi^i$ the $i$-th coordinate of $\xi$. Recall that $\Ha$ has drift $\frac{b^2}{2}$ and jump process $(\xi^1(t))$, and define $F(t):=\frac{b^2}{2}t+\sum_{s<t}\xi^1(s)$ for all $t\geq0$. Then 
\upshape$$\pi(\d z)=m\left(F(T)\in\d z,\ \{\xi^3(s)=0\ \forall s<T\right\}),$$
\itshape where $T:=\inf\{t\geq0,\ \xi^2(t)>x+F(t)\}$.
\end{remarque}

\par\bigskip
We turn our attention to the proof of Theorem \ref{th_cv_loi_des_sauts}, which will mainly rely on the following proposition.

\begin{proposition} \label{prop_cv_TL}
 The Laplace transform \upshape
$\E(e^{-r\Han(\en-)},\ L(\Txn)<\en)$
\itshape converges to \upshape
$$\E(e^{-r\Ha(\e-)},\ L(\Tx)<\e)=\frac{\lambda}{\lambda+\kill+\psi^*(r)}\int_{[0,\infty)} e^{-ar}\gamma^x(a,0,r)\pi(\d a),$$
\itshape where $\psi^*$ is the Laplace exponent of $H^*$ defined in Theorem \ref{th_cv_H*}, $\pi$ and $g^x$ are defined in Theorem \ref{th_cv_loi_des_sauts} above, and \upshape$\gamma^x(a,q,r):=\int_{[0,\infty)} e^{-rv}g^x(a,\{q\},\d v)$.
\end{proposition}

To prove the theorem and proposition above, we will need the following lemmas.

\begin{lemme} \label{lemme_cv_pi}
 Define for $a,h,t\in\R_+$: 
\upshape $$\pi_n(\d a):=\P(\bar \tZ_n(\Txn)\in\d a,\ L_n(\Txn)\leq\en).$$
\itshape Then $(\pi_n)$ converges weakly towards the measure $\pi$ defined in Theorem \ref{th_cv_loi_des_sauts}.
\end{lemme}

\begin{demo}
To prove the lemma we prove that $(\bar \tZ_n(\Txn),L_n(\Txn),\en)$ converges in distribution towards $(Z(\Tx),L(\Tx),\e)$. From Theorem \ref{th_appendix}, we know that the triplet $(\tZ_n,L_n,\Hcn)$ converges in distribution towards $(Z,L,\Hc)$. Using the Skorokhod representation theorem, there exists a sequence $(\mathcal \tZ_n,\mathcal L_n,\mathcalHcn)$  converging almost surely towards $(Z,L,\Hc)$, and such that $(\mathcal \tZ_n,\mathcal L_n,\mathcalHcn)$ and $(\tZ_n,L_n,\Hcn)$ are equal in law. We will use the notation $\mathcal T_n^{-x}$ for the first entrance time of $\mathcal \tZ_n$ in $\{-x\}$, and $\bar{\mathcal Z}_n(t)=\underset{[0,t]}\sup \mathcal \tZ_n$. \\

Thanks to Proposition \ref{prop_cv}.\eqref{prop_cv_Tx}, we know that as $n\to\infty$, $\mathcal T^{-x}_n\to \Tx$ a.s. Then note that $\mathcal T^{-x}_n$ is a continuity time for $\bar{\mathcal Z}_n$ and $\mathcal L_n$, since they are a.s. constant in a neighbourhood of $\mathcal T_n^{-x}$, and hence we get from Proposition 2.1 (b.5) in \cite{JS} that $\mathcal \tZ_n(\mathcal T^{-x}_n)\to Z(\Tx)$ and $\mathcal L_n(\mathcal T^{-x}_n)\to L(\Tx)$ a.s. \\

We have $\mathcal E_n=T^1(\mathcalHcn)$ and $\e=T^1(\Hc)$, where $\mathcalHcn$ and $\Hc$ are Poisson processes satisfying $\mathcalHcn\underset\P{\xrightarrow{a.s.}}\Hc$. Here Proposition VI.2.11 in\cite{JS} cannot be applied, although $\mathcal E_n$ is a first entrance time. But with an analogous proof, and using the fact that $\mathcalHcn$ is a Poisson process, we easily show that $\mathcal E_n\to\e$ a.s.\\ 

So, we have obtained the a.s. convergence (and thus the convergence in probability) 
\[(\bar{\mathcal Z}_n(\mathcal T_n^x),\mathcal L_n(\mathcal T_n^x),\mathcal E_n)\underset\P{\xrightarrow{a.s.}}(\bar Z(\Tx),L(\Tx),\e)\]
which gives, together with the equality in law  
$(\bar{\mathcal Z}_n(\mathcal T_n^x),\mathcal L_n(\mathcal T_n^x),\mathcal E_n)\stackrel{(d)}{=}(\bar \tZ_n(\Txn),L_n(\Txn),\en),$ 
the joint convergence in distribution of $(\bar \tZ_n(\Txn),L_n(\Txn),\en)$ towards $(Z(\Tx),L(\Tx),\e)$.
\end{demo}

\begin{lemme}
\label{lemme_g}
For all $y>0$, $v>0$ and $q\in\{0,1\}$, \upshape
\begin{align*}
 &\P_{-y}(Z(T^{(0,\infty)})\in\d v,\ \Delta\Zb(T^{(0,\infty)})\in\d q)\\
&=\frac{b^2}2 (W'(y)-\eta W(y))\delta_0(\d v)\delta_0(\d q)+\int_{(0,\infty)} (e^{-\eta u}W(y)-W(y-u))\;\B{f(u+v)}(\d q)\;\L(u+\d v)\;\d u,
\end{align*}
 \itshape and for all $n\geq1$ \upshape
\begin{align*}
 \P_{-y}(\tZ_n(T_n^{(0,\infty)})\in\d v,\ &\Delta\Zbn(T_n^{(0,\infty)})\in\d q)\\
&=\int_{(0,\infty)} (e^{-\teta_n v}\tilde W_n(y)-\tilde W_n(y-u))\;\B{f_n(n(u+v))}(\d q)\;\tL_n(u+\d v)\;\d u.
\end{align*}
\end{lemme}

The first quantity corresponds in fact to $g^x(y-x,\d q,\d v)$ introduced in the statement of Theorem \ref{th_cv_loi_des_sauts}. To keep consistency in the notation, we will then set
$$g_n^x(a,\d q,\d v):=\int_{(0,\infty)} (e^{-\teta_n v}\tilde W_n(x+a)-\tilde W_n(x+a-u))\B{f_n(n(u+v))}(\d q)\tL_n(u+\d v)\d u,$$
which corresponds to the second formula of Lemma \ref{lemme_g}.\\

\begin{demo}
We first write
\begin{align*}
 &\P_{-y}(\tZ_n(T_n^{(0,\infty)})\in\d v,\ \Delta\Zbn(T_n^{(0,\infty)})\in\d q)\\
&=\int_{[0,\infty)} \P_{-y}(\tZ_n(T_n^{(0,\infty)})\in\d v,\ \tZ_n(T_n^{(0,\infty)}-)\in\d u) \B{f_n(n(u+v))}(\d q).
\end{align*}
and similarly for $Z$.\\
Now according to \cite{K} (see consequence of (8.29)), we have for all $u>0$, $v>0$ :
$$\P_{-y}(Z(T^{(0,\infty)})\in\d v,\ Z(T^{(0,\infty)}-)\in\d u)=(e^{-\eta u}W(y)-W(y-u))\L(u+\d v)\d u,$$
and similarly
\begin{align*}
\P_{-y}(\tZ_n(T_n^{(0,\infty)})\in\d v,\ \tZ_n(T_n^{(0,\infty)}-)\in\d u)=(e^{-\teta_n v}\tilde W_n(y)-\tilde W_n(y-u))\tL_n(u+\d v)\d u. 
\end{align*}
Moreover, \cite[Exercise 8.6]{K} provides a formula for the probability of creeping over $0$ starting at $-y<0$ for a spectrally positive Lévy process, that is, the probability that the process is equal to $0$ at $T^{(0,\infty)}$ under $\P_{-y}$. In particular this probability is zero if the process has no Gaussian component, so that at rank $n$ we have
\begin{align*}
 &\P_{-y}(\tZ_n(T_n^{(0,\infty)})\in\d v,\ \Delta\Zbn(T_n^{(0,\infty)})\in \d q)\\
&=\int_{(0,\infty)} (e^{-\teta_n v}\tilde W_n(y)-\tilde W_n(y-u))\tL_n(u+\d v)\,\d u\;\B{f_n(n(u+v))}(\d q).
\end{align*}
On the other hand, as far as $Z$ is concerned, its Gaussian coefficient $\frac{b^2}2$ might be positive, and since $f(0)=0$, the formula of Exercise 8.6 in \cite{K} :
$$\P_{-y}(Z(T^{(0,\infty)})=0)=\frac{b^2}2(W'(y)-\eta W(y))$$
implies
\begin{align*}
 &\P_{-y}(Z(T^{(0,\infty)})\in\d v,\ \Delta\Zb(T^{(0,\infty)})\in \d q)\\
&=\frac{b^2}2 (W'(y)-\eta W(y)) \delta_0(\d v)\delta_0(\d q)+\int_{(0,\infty)} (e^{-\eta u}W(y)-W(y-u))\L(u+\d v)\,\d u\;\B{f(u+v)}(\d q),
\end{align*}
which ends the proof.
\end{demo}

\begin{lemme} \label{lemme_cv_gamma}
For all $n\geq1$, $y\in\R_+$, $r\in\R_+$, $q\in\{0,1\}$, define 
\upshape $$\gamma^x_n(y,q,r):=\int_{(0,\infty)} e^{-rv}g_n^x(y,\{q\},\d v).$$
\itshape Then the Laplace transform $\gamma_n^x(y,q,r)$ converges towards $\gamma^x(y,q,r)$ (defined in Proposition \ref{prop_cv_TL}) as $n\to\infty$, and the convergence is uniform w.r.t. $y$ on every compact set of $\R_+$.
\end{lemme}

\begin{demo}
Fix $y\in\R_+$, $r\in\R_+$, $q\in\{0,1\}$. Using the expression of $g^x$ given by formula (\ref{formula_g}), we have :
$$ \gamma^x_n(y,q,r) =\E_{-(x+y)}\big(e^{-r\tZ_n(T_n^{(0,\infty)})},\ \Delta \Zbn(T_n^{(0,\infty)})=q\big) ,$$
which we also can reexpress as :
\begin{align*}
 \gamma^x_n&(y,q,r)  =\E\big(e^{-r(\Han(L_n(T_n^{(x+y,\infty)}))-(x+y))},\ \Delta \Hcn(L_n(T_n^{(x+y,\infty)}))=q\big) \\
		   & =\int_{(x+y,\infty)\times(0,\infty)} \B{f_n(nu)}(\d q)\, e^{-r(v-(x+y))}\, \P(\Han(L_n(T_n^{(x+y,\infty)}))\in\d v,\Delta \tZ_n(T_n^{(x+y,\infty)})\in\d u).
\end{align*}
In the same vein, we have 
$$\gamma^x(y,q,r) =\int_{(x+y,\infty)\times[0,\infty)} \B{f(u)}(\d q)\, e^{-r(v-(x+y))}\, \P(\Ha(L(T^{(x+y,\infty)}))\in\d v,\Delta Z(T^{(x+y,\infty)})\in\d u).$$

To start with, we prove that the measures $\P(\Han(L_n(T_n^{(x+y,\infty)}))\in\d v,\ \Delta \tZ_n(T_n^{(x+y,\infty)})\in\d u)$ converge weakly towards $\P(\Ha(L(T^{(x+y,\infty)}))\in\d v,\ \Delta Z(T^{(x+y,\infty)})\in\d u)$. First recall that thanks to Theorem \ref{th_appendix} we have the convergence in distribution of $(\Han,\tZ_n)$ towards $(\Ha,Z)$. With probability one, we have that $T^{(x+y,\infty)}$ is either a continuity point of $Z$ a.s., or it satisfies $Z(T^{(x+y,\infty)}-)<x+y<Z(T^{(x+y,\infty)})$ and $\Ha(L(T^{(x+y,\infty)})-)<x+y<\Ha(L(T^{(x+y,\infty)}))$ a.s. Note furthermore that the first entrance time of $\Ha$ in $(x+y,\infty)$ is equal to $L(T^{(x+y,\infty)})$ a.s. Then we can easily adapt the proof of Proposition VI.2.12 in \cite{JS} to get that 
$$(\Han(L_n(T_n^{(x+y,\infty)})),\ \Delta \tZ_n(T_n^{(x+y,\infty)})) \Rightarrow (\Ha(L(T^{(x+y,\infty)})),\ \Delta Z(T^{(x+y,\infty)})).$$
On the other hand, under \HypMutConst\ as well as under \HypMut, we have the uniform convergence of $(u,v)\mapsto \B{f_n(nu)}\{q\}\, e^{-rv}$ to $(u,v)\mapsto \B{f(u)}\{q\}\, e^{-rv}$ on every compact set of $(0,\infty)\times(x+y,\infty)$. Then from an appeal to Lemma \ref{lemme_Cb} we get the convergence of $\gamma^x_n(y,q,r)$ towards  $\gamma^x(y,q,r)$ for all fixed $y,r\geq0$ and $q\in\{0,1\}$.\\

It remains to prove the uniform convergence of $\gamma_n^x$ w.r.t. the first variable, $y$, on every compact set of $\R_+$.\\
Take $r\geq0$ and $q\in\{0,1\}$. For all $y\geq0$, $\tilde W_n(x+y)$ is positive, and we set $\tilde \gamma_n^x(y,q,r)=\gamma_n^x(y,q,r) / \tilde W_n(x+y)$. Observe that $y\mapsto\tilde \gamma_n^x(y,q,r)$ is decreasing on $\R_+$ : Indeed,  
it can be shown with elementary calculations that
$$\tilde\gamma_n^x(y,q,r)=\int_{(0,\infty)} \tL_n(\d z) \B{f_n(nz)}(\{q\}) \int_0^z e^{-r(z-u)}\ e^{-\teta_n u} N_n(\inf\epsilon \leq -(x+y) \sachant -\epsilon(\zeta-)=u) \ \d u,$$
and since the mappings $y\mapsto N_n(\inf\epsilon \leq -(x+y) \sachant -\epsilon(\zeta-)=u)$ are clearly decreasing, we have the same property for $\tilde\gamma_n^x(\cdot,q,r)$. Next, recalling that the functions $\tilde W_n(x+\cdot)$ are strictly increasing and take positive values, we get that the functions $\tilde\gamma_n^x(\cdot,q,r)$ are decreasing, which leads to the uniform convergence of $\gamma_n^x(\cdot,q,r)$ to $\gamma^x(\cdot,q,r)$ on every compact set of $\R_+$.
\end{demo}

In the proof of Proposition \ref{prop_cv_TL}, we will make a frequent use of the Markov property, applied alternately to $Z$ (resp. $\tZ_n$) or to $(\Ha,\Hc)$ (resp. $(\Han,\Hcn)$), at different stopping times. We already know that $\Tx$, $\Linv(t)$ (resp. $\Txn$, $\Linv_n$) are $\F$- (resp. $\F_n$-) stopping times. We introduce here three other stopping times which we will need later.\\

First we define the processes $\Zcn,\Zc$ as follows : for all $t\geq0$,
$$\Zcn(t):=\Hcn(L_n(t)-)\ \text{ and }\ \Zc(t)=\Hc(L(t)-).$$
The process $\Zcn$ is a counting process which jumps every time a mutation occurs at a record time of $\tZ_n$ : it can be seen as the matching process of $\Hcn$ in the real time scale (in opposition to the local time scale) - and similarly for $\Zc$. We then have the identity $\Hcn=\Zcn\circ\Linv_n$.\\
We enlarge the initially considered filtrations and set for all $t\geq0$ :
$$\Fc_t:=\sigma(Z(s),\Zc(s),\, s\leq t)$$
and
$$\Fc_{n,t}:=\sigma(\tZ_n(s),\Zcn(s),L_n(s)\, s\leq t).$$

We denote by $\Fc$ (resp. $\Fc_n$) the filtration $(\Fc_t)_{t\geq0}$ (resp. $(\Fc_{n,t})_{t\geq0}$). Finally, the notations $\Fc_{\Linv}$, $\Fc_{n,\Linv_n}$ will respectively stand for the filtrations $(\Fc_{\Linv(t)})_{t\geq0}$ and $(\Fc_{n,\Linv_n(t)})_{t\geq0}$. Note that $(Z,\Zc)$ is not a Markov process in the filtration $\Fc$.\\

\begin{proposition}\label{prop_StoppingTimes}
 \begin{enumerate}[\upshape(i)]
  \item \upshape $\e$ \itshape (resp. \upshape $\en$) \itshape is a stopping time w.r.t. the filtration $\Fc_{\Linv}$ (resp. $\Fc_{n,\Linv_n}$).
  \item \upshape $\Linv(\e-)$ \itshape (resp. \upshape $\Linv_n(\en-)$) \itshape is a stopping time w.r.t. the filtration $\Fc$ (resp. $\Fc_n$).
  \item $L(\Tx)$ (resp. $L_n(\Txn)$) \itshape is a stopping time w.r.t. the filtration $\Fc_{\Linv}$ (resp. $\Fc_{n,\Linv_n}$).
 \end{enumerate}
\end{proposition}

\begin{figure}[!h]
\begin{center}
  \includegraphics[width=15.5cm]{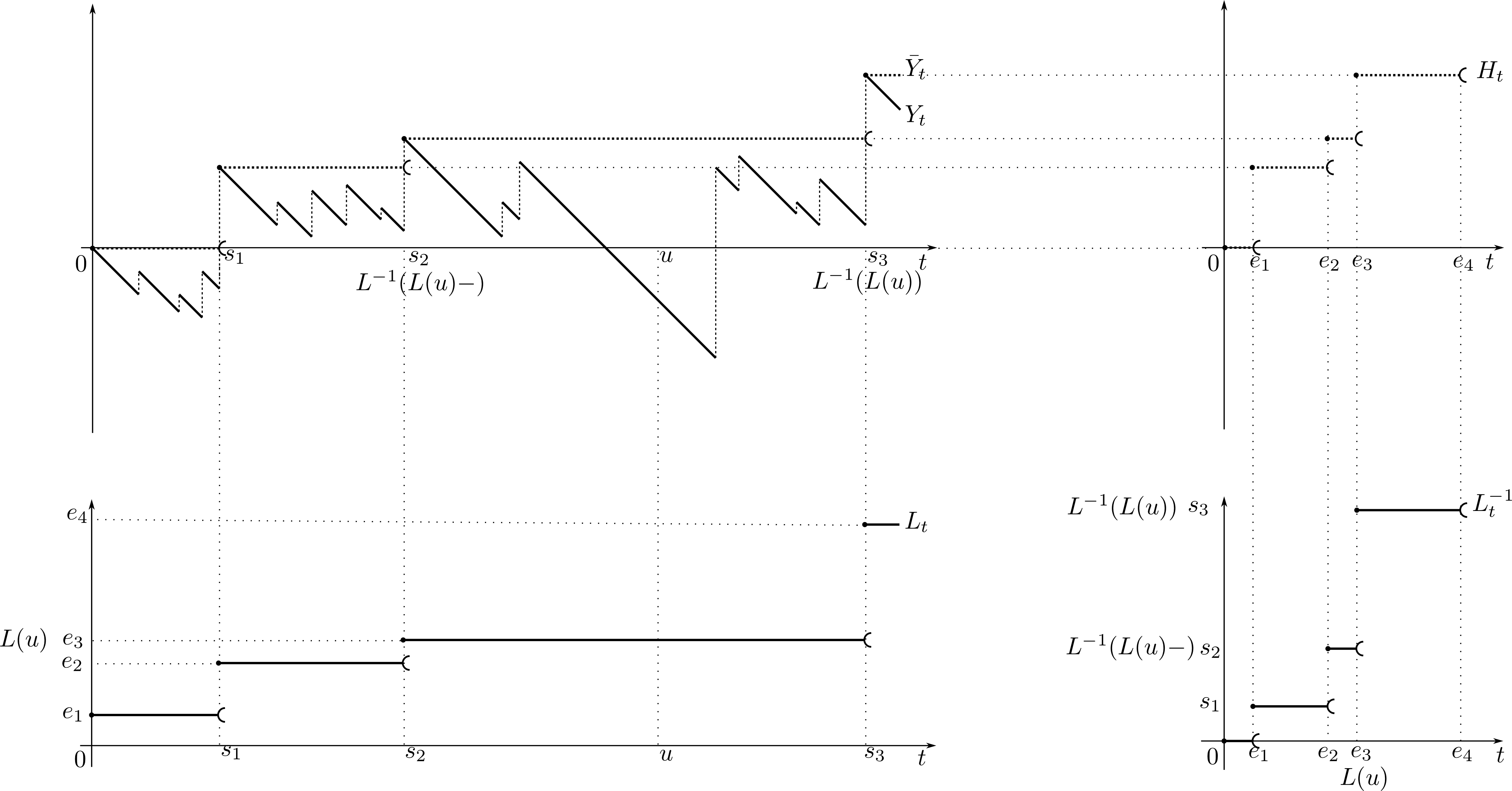}
 \caption{An example of the paths of $Z$, its local time at the supremum and its ladder process in the finite variation case.}
\label{figure_LP}
\end{center}
 \end{figure}

\begin{demo}
\begin{enumerate}[(i)]
 \item $\e$ is the first entrance time of the bivariate subordinator $(\Ha,\Hc)$ into $\R_+\times\R_+^*$. Thus $\e$ is a stopping time w.r.t. the natural filtration associated to $(\Ha,\Hc)$, and then w.r.t. $\Fc_{\Linv}$.
 \item We have for all $t\geq0$ :
$$\{\Linv(\e -)\leq t\}=\bigcap_{u\in\mathbb Q\cap\R_+^*}\big(\{\Linv(u)\leq t\}\cap\{u<\e\}\big),$$
now $\Linv(t)$ is a $\F$- stopping time, thus $\{\Linv(u)\leq t\}\in\F_t$, and $\e$ is a $\Fc_{\Linv}$-stopping time, thus $\{\e>u\}\in\Fc_{\Linv(u)}$. Consequently $\{\Linv(u)\leq t\}\cap\{u<\e\}$ belongs to $\Fc_t$ for all $u>0$, and so is $\{\Linv(\e -)\leq t\}$.
 \item For all $t\geq0$, we want to prove that $\{L(\Tx)\leq t\}=\{\Tx\leq \Linv(t)\}$ a.s. For a clearer view of what follows, see Figure \ref{figure_LP}.\\
Fix $u\geq0$. On the one hand, since $u\leq \Linv(L(u))$ and $\Linv$ is increasing, $L(u)\leq t$ implies $u\leq\Linv(t)$. On the other hand, in the infinite variation case, the function $L\circ\Linv$ is the identity function, and hence $u\leq \Linv(t)$ implies $L(u)\leq t$. In the finite variation case, the definition of $\Linv$ implies that if $u<\Linv(t)$, then $L(u)\leq t$. Now the event $\{\exists t\geq0,\ \Tx=\Linv(t)\}$ is negligible, thus $\{\Tx\leq\Linv(t)\}=\{\Tx<\Linv(t)\}$ a.s. \\
We conclude from what precedes that the events $\{L(\Tx)\leq t\}$ and $\{\Tx\leq \Linv(t)\}$ are identical a.s., and since $\Tx$ is a stopping time w.r.t. the filtration $\F$, this implies that $L(\Tx)$ is a stopping time w.r.t. the filtration $\Fc_{\Linv}$.
\end{enumerate}
Remark that the three proofs above work in the infinite variation case as well as in the finite variation case, so that the conclusions are also true for $\en$, $\Linv_n(\en-)$, $L_n(\Txn)$, $n\geq1$.

\end{demo}

\par\bigskip
\begin{demopr}{\ref{prop_cv_TL}}
We begin with the computation of the probability measure $\P(\Ha(\e-)\in\d z,\ L(\Tx)<\e)$ :  The calculation below is done for the limiting process. However we pay attention to the fact that the arguments are still valid in the finite variation case, so that the same calculation remains true for $\P(\Han(\en-)\in\d z,\ L_n(\Txn)<\en)$. \\

Noting that $L(\Tx)<\e$ coincides with $\Hc(L(\Tx))=0$ a.s., and then applying the Markov property to the process $(\Ha,\Hc)$ at the $\Fc_{L^{-1}}$-stopping time $L(\Tx)$, we have
\begin{align*}
&\P(\Ha(\e-)\in\d z,\ L(\Tx)<\e) \\
&=\int_{[0,z)} \int_{(a,z)} \P_b(\Ha(\e-)\in\d z)\, \P(\Ha(L(\Tx)-)\in\d a,\ \Ha(L(\Tx))\in\d b,\ \Hc(L(\Tx))=0)
\end{align*}
Using the notation $d(\Tx):=\Linv(L(\Tx))$, recall that 
$$\Ha(L(\Tx)-)=\bar Z(\Tx) \text{ and } \Ha(L(\Tx))=\bar Z(d(\Tx)).$$ 
Furthermore, with probability one 
$$\{\Hc(L(\Tx))=0\}=\{\Zc(\Tx)=0\}\cap\{\Delta\Zc(d(\Tx))=0\}.$$ 
Conditional on $\bar Z(\Tx)$, the random variable $\Delta\Zc(d(\Tx))$ is independent from $\Fc_{\Tx}$, and has the law of $\Delta\Zb(T^{(0,\infty)})$ under $\P_{-x-\bar Z(\Tx)}$ (Note that $T^{(0,\infty)}$ is a.s. necessarily a record time for $Z$ under $\P_{-x-\bar Z(\Tx)}$). We then use the Markov property again, applied to the process $Z$ at the $\F$-stopping time $\Tx$ :
\begin{multline*}
\P(\Ha(L(\Tx)-)\in \d a,\ \Ha(L(\Tx))\in\d b,\ \Hc(L(\Tx))=0)\\
\qquad\quad=\P(\bar Z(\Tx)\in\d a,\ \Zc(\Tx)=0)\; \P_{-(a+x)}(Z(T^{(0,\infty)})\in\d b-a,\ \Delta\Zb(T^{(0,\infty)})=0)\\
=\P(\bar Z(\Tx)\in\d a,\ L(\Tx)\leq\e)\; \P_{-(a+x)}(Z(T^{(0,\infty)})\in\d b-a,\ \Delta\Zb(T^{(0,\infty)})=0),
\end{multline*}
and finally, using the notation introduced in the statement of Theorem \ref{th_cv_loi_des_sauts} and Lemma \ref{lemme_g}, this gives
\begin{multline*}
\P(\Ha(\e-)\in\d z,\ \Delta\Ha(\e)\in\d y,\ L(\Tx)<\e) \\
=\P(\Delta\Ha(\e)\in\d y)\;\int_{[0,z)}\int_{(a,z)} \P_b(\Ha(\e-)\in\d z)\;\pi(\d a)\; g^x(a,\zero,\d b-a).
\end{multline*}

When $Z$ does not drift to $-\infty$, by definition of $\e$ as first entrance time and thanks to Proposition 0.5.2 in \cite{B}, we have 
$$\P_b(\Ha(\e-)\in\d z)=\P_b(H^*(\alpha)\in\d z),$$
where $H^*$ is the subordinator defined in Theorem \ref{th_cv_H*}, and  $\alpha$ is an independent exponential random variable with parameter $\lambda$. 

In the same way, we treat the case $Z$ drifts to $-\infty$ appealing to \cite[Prop. 0.5.2]{B} and to Theorem \ref{th_cv_H} : set $\varrho:=\e\wedge\mathcal K$, where $\mathcal K$ is an independent exponential variable with parameter $\kill$. Then $\varrho$ follows an exponential distribution with parameter $\lambda+\kill$, and we have 
$$\P_b(\Ha(\e-)\in\d z)=\frac{\lambda}{\lambda+\kill}\P_b(\Ha(\varrho-)\in\d z\sachant \varrho=\e)=\frac{\lambda}{\lambda+\kill}\P_b(H^*(\alpha')\in\d z),$$
where $\alpha'$ is an independent exponential random variable with parameter $\lambda+\kill$. \\

By definition of $U_*^{(\cdot)}$ we have then in both cases $\P_b(\Ha(\e-)\in\d z)=\lambda U_*^{(\lambda+\kill)}(\d z-b)$ (recall that we set $\kill=0$ if $Z$ does not drift to $-\infty$). As a consequence, 
\begin{align*}
\P(\Ha(\e-)\in\d z,\ L(\Tx)<\e) = \lambda\;\int_{[0,z)}\int_{(a,z)} U_*^{(\lambda+\kill)}(\d z-b)\;\pi(\d a)\; g^x(a,\zero,\d b-a).
\end{align*}

Hence we get for the Laplace transform :
\begin{align*}
\E(e^{-r\Ha(\e-)},\ L(\Tx)<\e) = \int_{[0,\infty)}\pi(\d a) \int_{[a,\infty)} g^x(a,\zero,\d b-a) \int_{[b,\infty)} e^{-rz}\lambda U_*^{(\lambda+\kill)}(\d z-b).
\end{align*}
From the definition of $U_*^{(\lambda+\kill)}$ we have for all $r\geq0$, $\int_{(0,\infty)} e^{-rz} U_*^{(\lambda+\kill)}(\d z)=(\lambda+\kill+\psi^*(r))^{-1}$, which leads to
\begin{align*}
\E(e^{-r\Ha(\e-)},\ L(\Tx)<\e)& = \frac{\lambda}{\lambda+\kill+\psi^*(r)}\int_{[0,\infty)}\pi(\d a) \int_{[a,\infty)} g^x(a,\zero,\d b-a) e^{-br}\\
& = \frac{\lambda}{\lambda+\kill+\psi^*(r)}\int_{[0,\infty)}\pi(\d a) \gamma^x(a,0,r) e^{-ar},
\end{align*}
and as announced, we have a similar formula at rank $n$ :
\begin{align*}
\E(e^{-r\Han(\en-)},\ L_n(\Txn)<\en)= \frac{\lambda_n}{\lambda_n+\kill_n+\psi_n^*(r)}\int_{[0,\infty)}\pi_n(\d a) \gamma_n^x(a,0,r) e^{-ar}.
\end{align*}

Now as $n\to\infty$, thanks to \cite[Prop. 4.9.(i)]{LPWM1} $\lambda_n=\mu_n(\R_+^*,\un)$ converges to $\lambda$, and thanks to Theorem \ref{th_cv_H*} $\psi^*_n$ converges to $\psi^*$. According to the proof of \ref{th_cv_H} in \cite{LPWM1}, we also have $\kill_n\to\kill$. As for the integral, thanks to Lemma \ref{lemme_cv_pi} and Lemma \ref{lemme_cv_gamma} we can apply Lemma \ref{lemme_Cb}, and hence we have proved that as $n\to\infty$,
$$\E(e^{-r\Han(\en-)},\ L_n(\Txn)<\en) \to \E(e^{-r\Ha(\e-)},\ L(\Tx)<\e)$$ for all $r\geq0$. This finishes the proof.
\end{demopr}

\par\bigskip
Finally, before we prove the theorem, we need the following technical lemma :
\begin{lemme} \label{lemme_events}
 The event \upshape $\{L(\Tx)<\e\}$ \itshape (resp. \upshape$\{L_n(\Txn)<\en\}$) \itshape belongs to \upshape  $\Fc_{\Linv(\e-)}$\itshape (resp. \upshape$\Fc_{n,\Linv_n(\en-)}$\itshape).
\end{lemme}
\begin{demo}
We first want to prove that $\{L(\Tx)<\e\}=\{\Linv(L(\Tx))<\Linv(\e)\}$ a.s. (and the equivalent equality at rank $n$). \\

As far as the limiting process is concerned, we are in the infinite variation case : The process $\Linv$ is a.s. continuous and strictly increasing, so that $\{L(\Tx)<\e\}=\{\Linv(L(\Tx))<\Linv(\e)\}$ a.s. In fact these two events still coincide a.s. in the finite variation case, although $\Linv$ is not strictly increasing : Indeed, 
$\Linv$ is injective on the set of all jumping times of $\Ha$ ; now $L(\Tx)$ and $\e$ are a.s. two jumping times of $\Ha$, hence $\Linv(L(\Tx))=\Linv(\e)$ implies $L(\Tx)=\e$, and the claim is proved. \\

We now prove that $\Linv(L(\Tx))<\Linv(\e) \Leftrightarrow \Tx\leq\Linv(\e-)$ a.s.\\
On the one hand, $\Linv(L(\Tx))$, $\Linv(\e)$ and $\Linv(\e-)$ belong to the zero set of $\bar Z-Z$, and $\Linv(\e)$ and $\Linv(\e-)$ are two consecutive (possibly equal) zeros of $\bar Z-Z$. Thus $\Linv(L(\Tx))<\Linv(\e)$ implies $\Linv(L(\Tx))\leq\Linv(\e-)$, and since $\Tx\leq\Linv(L(\Tx))$ a.s., this ensures that with probability one $\{\Linv(L(\Tx))<\Linv(\e)\}\subset\{\Tx\leq\Linv(\e-)\}$. 

On the other hand, assume that $\Tx\leq\Linv(\e-)$. The event $\{\Tx=\Linv(t), \text{ for some } t\geq0\}$ is negligible and thus by definition of $\Linv(\e-)$, there exists $u<\e$ such that $\Tx<\Linv(u)$ a.s. This ensures that $\Linv(L(\Tx))=\inf\{\Linv(u),\ \Linv(u)>\Tx\}<\Linv(\e)$ a.s., and then $\{\Tx\leq\Linv(\e-)\}\subset\{\Linv(L(\Tx))<\Linv(\e)\}$\\

So, we have proved that almost surely $\{L(\Tx)<\e\}=\{\Tx\leq\Linv(\e-)\}$ a.s. We conclude using the fact that $\Tx$ is a $\Fc$-stopping time. The proof above remains true in the finite variation case, so that the result is also valid at rank $n$.
\end{demo}

\begin{demoth}{\ref{th_cv_loi_des_sauts}}
To begin with, we prove formulas \eqref{formula_cv} and \eqref{formula_cv_Const}. As in the proof above, we do the calculation and reasoning for the limiting process $Z$, and we add some remarks when needed so that it remains valid at rank $n$.\\

First note that thanks to the Markov property applied to $(\Ha,\Hc)$ at the $\Fc_{\Linv}$-stopping time $\e$, and since the event $\{\Linv(\e)<\Tx\}=\{\Tx\leq\Linv(\e)\}^c$ (where $A^c$ denotes the complementary event of $A$) belongs to $\Fc_{\Linv(\e)}$, we have
\begin{align*}
& \P(\Ha(\e-)\in\d z,\ \Delta\Ha(\e)\in\d y,\ \Linv(\e)<\Tx<T^{(\tau-x,\infty)})\\
&=\P(\Ha(\e-)\in\d z,\ \Delta\Ha(\e)\in\d y,\ \Linv(\e)<\Tx)\P_{z+y}(\Tx<T^{(\tau-x,\infty)}),
\end{align*}
where $\P_{z+y}(\Tx<T^{(\tau-x,\infty)})=W(\tau-x-z-y)/W(\tau)$. Now we have :
\begin{align*}
 \P(\Ha(\e-)\in\d z,\ & \Delta\Ha(\e)\in\d y,\ \Linv(\e)<\Tx)\\
=\P(\Ha(\e-)\in\d z,\ &\Delta\Ha(\e)\in\d y)-\P(\Ha(\e-)\in\d z,\ \Delta\Ha(\e)\in\d y,L(\Tx)\leq\e)\\
=\P(\Ha(\e-)\in\d z,\ &\Delta\Ha(\e)\in\d y) -\P(\Ha(\e-)\in\d z,\ \Delta\Ha(\e)\in\d y,L(\Tx)<\e),\\
&-\P(\Ha(\e-)\in\d z,\ \Delta\Ha(\e)\in\d y,\ L(\Tx)=\e),
\end{align*}
where in the last equality we distinguished the case where the first mutation, in the time scale of $Z$, occurs at the end of the excursion interval containing $\Tx$, or later.  Recall  that despite the fact that $\Linv$ shall not be strictly increasing (finite variation case), we always have $L(\Tx)<\e\Leftrightarrow \Linv(L(\Tx))<\Linv(\e)$ (see proof of Lemma \ref{lemme_events}).\\

As in the proof of Proposition \ref{prop_cv_TL}, applying Proposition 0.5.2 in \cite{B}, we get for the first term in the sum :
$$\P(\Ha(\e-)\in\d z,\Delta\Ha(\e)\in\d y)=\lambda U_*^{(\lambda+\kill)}(\d z)\P(\Delta \Ha(\e)\in\d y),$$

Then we compute the second term in the sum : Lemma \ref{lemme_events} ensures that $\{L(\Tx)<\e\}\in\Fc_{\Linv(\e-)}$ a.s., thus by Markov property applied to $Z$ at the $\Fc$-stopping time $\Linv(\e-)$ we have 
\begin{align*}
 \P(\Ha(\e-)\in\d z,\,\Delta\Ha(\e)\in\d y,\,L(\Tx)<\e) =\P(\Ha(\e-)\in\d z,\,L(\Tx)<\e)\P(\Delta\Ha(\e)\in\d y),
\end{align*}
and thanks to the calculation made in the proof of Proposition \ref{prop_cv_TL}, we get 
\begin{align*}
\P(\Ha(\e-)\in&\d z,\ \Delta\Ha(\e)\in\d y,\ L(\Tx)<\e) =\\
& \lambda\P(\Delta \Ha(\e)\in\d y)\;\int_{[0,z)}\int_{(a,z)} U_*^{(\lambda+\kill)}(\d z-b)\;\pi(\d a)\; g^x(a,\zero,\d b-a).
\end{align*}

Finally note that $\lambda\P(\Delta \Ha(\e)\in\d y)=\theta\delta_0(\d y)$ under Assumption \HypMutConst, $\lambda\P(\Delta \Ha(\e)\in\d y)=\mu(\d y,\un)+\rho \delta_0(\d y)$ under Assumption \HypMut, and $\lambda_n\P(\Delta \Han(\en)\in\d y)=\mu_n(\d y,\un)$ in both cases. \\

It remains to compute the third term in the sum. With an application of the Markov property to $Z$ at $\Tx$ as in the proof of Proposition \ref{prop_cv_TL}, we get 
$$\P(\Ha(\e-)\in\d z,\ \Delta\Ha(\e)\in\d y,\ L(\Tx)=\e)=\pi(\d z)g^x(z,\un,\d y),$$
which vanishes in case \HypMutConst\ according to Remark \ref{remark_gx}.\\

Thereby we have established formulas \eqref{formula_cv} and \eqref{formula_cv_Const}, and these formulas remain true at rank $n$ (considering in case \HypMut\ that the coefficient $\rho$ is zero in the finite variation case). Since the expression of $g^x$ given by formula (\ref{formula_g}) in the statement of the theorem has been established in Lemma \ref{lemme_g}, proving the claimed convergence will end the proof.\\

From the calculation above we have at rank $n$ (in both cases \HypMutConst\ and \HypMut) :
\begin{multline}\label{formula_cv_n}
 \P(\Han(\en-)\in\d z,\ \Delta\Han(\en)\in\d y,\ \Linv_n(\en)<\Txn<T_n^{(\tau-x,\infty)}) \\
=\frac{\tilde W_n(\tau-z-y)}{\tilde W_n(\tau)} \Bigg\{\mu_n(\d y,\un)\Big[ U_{n,*}^{(\lambda_n+\kill_n)}(\d z)-\P(\Han(\en-)\in\d z,\ L_n(\Txn)<\en)\Big]\\
-\pi_n(\d z)g_n^x(z,\un,\d y)\Bigg\},
\end{multline}
where $\pi_n$, $g_n^x$ have been defined respectively in Lemmas \ref{lemme_cv_pi} and \ref{lemme_cv_gamma}, and $U_{n,*}^{(\lambda_n+\kill_n)}$ denotes the $(\lambda_n+\kill_n)$-resolvent measure of $H_n^*$ (defined in Theorem \ref{th_cv_H*}). Now as $n\to\infty$, for all $z\geq0$, $y>0$,

\begin{itemize}
 \item From Proposition \ref{prop_cv}.(\ref{prop_cv_W}), we know that $\tilde W_n(\tau-x-z-y)/\tilde W_n(\tau)$ converges to $W(\tau-x-z-y)/W(\tau)$.
 \item From \cite[Prop. 4.9.(i)]{LPWM1}, $\lambda_n=\mu_n(\R_+^*,\un)$ converges to $\lambda$, and $\mu_n(\d y,\un)$ converges weakly to $\theta\delta_0(\d y)$ (resp. $\mu(\d y,\un)+\rho\delta_0(\d y)$) in case \HypMutConst\ (resp. \HypMut).
 \item The Laplace transform of the measure $U_{n,*}^{(\lambda_n+\kill_n)}$ (resp. $U_{*}^{(\lambda+\kill)}$) is given by $(\lambda_n+\kill_n+\psi_n^*(\cdot))^{-1}$ (resp. $(\lambda+\kill+\psi^*(\cdot))^{-1}$), hence the measure $U_{n,*}^{(\lambda_n+\kill_n)}$ converges weakly towards $U_{*}^{(\lambda+\kill)}$ using Theorem \ref{th_cv_H*} and the fact that $\kill_n\to\kill$ (see proof of Theorem \ref{th_cv_H} in \cite{LPWM1}).
 \item The weak convergence of the probability measure $\P(\Han(\en)\in\d z,\ L_n(\Txn)<\en)$ has been proved via the convergence of its Laplace transform in Proposition \ref{prop_cv_TL}.
 \item Finally, the weak convergence of $\pi_n(\d z)g_n^x(z,\un,\d y)$ to $\pi(\d z)g^x(z,\un,\d y)$ is straightforward from Lemmas \ref{lemme_cv_pi} and \ref{lemme_cv_gamma}.
\end{itemize}
As a conclusion we have proved the weak convergence under Assumption \HypMut\ (resp. \HypMutConst) of (\ref{formula_cv_n}) to (\ref{formula_cv}) (resp. \eqref{formula_cv_Const}) .
\end{demoth}

\subsubsection{Convergence in distribution of $(M\neps)_n$ towards $M_\eps$} \label{sec_Proof_cvMeps}

The aim of this last subsection is to prove Proposition \ref{prop_cv_M}, appealing to Theorem 1 in \cite{Karr}. The four lemmas below ensure that the conditions needed to apply this theorem are fulfilled : First in Lemma \ref{lemme_cv4b_Karr}, we make use of Theorem \ref{th_cv_loi_des_sauts} to obtain a slightly more precise result about the convergence of the transition measures $\num$ and $\nud$.  Second, we give in Lemma \ref{lemme_formula_nui} explicit expressions for $\nuin$ and $\nui$, which allow us to prove in Lemma \ref{lemme_cv_nui} the weak convergence of $\nuin$ towards $\nui$. Finally we deduce from this the convergence in distribution of $M\neps(0)$ towards $M_\eps(0)$. \\

\begin{lemme}\label{lemme_cv4b_Karr}
 Suppose $x_n\to x$ as $n\to\infty$, where $x_n$ and $x$ are positive real numbers. Then $\num_n(x_n,\cdot)$ (resp. $\nud_n(x_n,\cdot)$) converges weakly towards $\num(x,\cdot)$ (resp. $\nud(x,\cdot)$).
\end{lemme}

\begin{demo}
Let $A$ be a set in $\mathcal B([\eps,\tau)\times\{0,1\})$ satisfying $\num(x,\partial A)=0$ (where $\partial A$ denotes the boundary of $A$). First write 
$$|\num_n(x_n,A)-\num(x,A)|\leq |\num_n(x_n,A)-\num_n(x,A)|+|\num_n(x,A)-\num(x,A)|.$$
The second term in the right-hand side vanishes thanks to Theorem \ref{th_cv_loi_des_sauts}. Besides, we have
$$|\P_0(\Han(\e_n)\in A, \Linv_n(\e_n)<T_{n}^{-x_n} )-\P_0(\Ha_n(\e_n)\in A, \Linv_n(\e_n)<\Tx)| \leq \P_0(\Txn < T_{n}^{-x_n} ),$$
which vanishes as $n\to\infty$ thanks to the a.s. continuity of $x\mapsto\Txn$ on $\R_+$ under $\P_0$. Then, by definition of $\num_n$ (see \eqref{formula_def_num}), we get $|\num_n(x_n,A)-\num_n(x,A)|\to0$ as $n\to\infty$ (for the sake of simplicity, we omitted here the conditioning that appears in the definition of $\num_n$).

A similar reasoning holds for $\nud$ (the weak convergence of $\nud_n$ towards $\nud$ is a consequence of Lemma \ref{lemme_cv_pi}).
\end{demo}

\par\bigskip

\begin{lemme}\label{lemme_formula_nui}
For all $(u,q)\in[\eps,\tau)\times\{0,1\}$, we have \upshape
\begin{multline}\label{formula_nuin}
 \nuin(\d u,\d q)\\
=\left(\frac1{\tilde W_n(\eps)}-\frac1{\tilde W_n(\tau)}\right)\frac{\tilde W_n(\tau-u)}{\tilde W_n(\tau)}  \d u \int_{(u,\infty)} \tL_n(\d z)\B{f_n(nz)}(\d q) \bigg(1-\frac{\tilde W_n(\eps-(z-u))}{\tilde W_n(\eps)}\bigg),
\end{multline}
\begin{multline}\label{formula_nui}
 \nui(\d u,\d q)\\
=\frac1{p_\eps} \frac{W(\tau-u)}{ W(\tau)} \left[\frac{b^2}{2}\frac{W'(\eps)}{W(\eps)}\delta_{\eps}(\d u)\delta_0(\d q)+\d u\int_{(u,\infty)} \L(\d z)\B{f(z)}(\d q) \bigg(1-\frac{W(\eps-(z-u))}{W(\eps)}\bigg)\right],
\end{multline}
\end{lemme}

 \begin{lemme}\label{lemme_cv_nui}
The sequence of measures $(\nuin)$ converges weakly towards $\nui$.
\end{lemme}

\par\bigskip
For the sake of clarity we only prove the two lemmas for $\nuin(\,\cdot\,,\{0,1\})$ and $\nui(\,\cdot\,,\{0,1\})$.  \\

\begin{demole}{\ref{lemme_formula_nui}}
We begin by proving \eqref{formula_nuin}. Recall the following definition :
$$\nuin(\,\cdot\,,\{0,1\})=\frac{1}{p'\neps}\P_0(\Ups_n\in\,\cdot\,,\ T_n^{-\eps}<T_n^{(0,\infty)}<T_n^{-\tau}).$$
Applying the strong Markov property at $T_n^{-\eps}$, we get :
$$ \P_0 (\Ups_n-\eps\in\d u,\ T_n^{-\eps}<T_n^{(0,\infty)}<T_n^{-\tau})=\P_0(T_n^{-\eps}<T_n^{(0,\infty)})\P_{-\eps}(\Ups_n-\eps\in\d u,\ T_n^{(0,\infty)}<T_n^{-\tau}).$$
Now conditional on $T_n^{(0,\infty)}<T_n^{-\tau}$, $\Ups_n-\eps$ has the distribution under $N'_n(\cdot\sachant -\inf\epsilon<\tau-\eps,\ \sup\epsilon\geq\eps)$ of the undershoot of an excursion at its first entrance time in $(0,\infty)$. Thanks to Proposition 0.5.2(ii) in \cite{B}, we then have 
\begin{multline*}
 \P_0 (\Ups_n-\eps\in\d u,\ T_n^{-\eps}<T_n^{(0,\infty)}<T_n^{-\tau}) \\
=\P_0(T_n^{-\eps}<T_n^{(0,\infty)})\P_{-\eps}(T_n^{(0,\infty)}<T_n^{-\tau})\frac{N'_n(-\epsilon(\chi-)\in\d u,\ -\inf\epsilon<\tau-\eps,\ \sup\epsilon\geq \eps)}{N'_n(-\inf\epsilon<\tau-\eps, \sup\epsilon\geq \eps)},
\end{multline*}
Recall that for any $\epsilon\in\mathscr E'$, $\chi(\epsilon)$ denotes the first (and unique) entrance time of $\epsilon$ into $(0,\infty)$, which is a.s. finite on $\{-\inf\epsilon<\tau-\eps\}$.

The process $\tZ_n$ has finite variation, and it can be shown with elementary calculations that $$N'_n(-\inf\epsilon<\tau-\eps, \sup\epsilon\geq \eps)=\frac{\tilde W_n(0)}{\tilde W_n(\tau-\eps)}\left(\frac{\tilde W_n(\tau)}{\tilde W(\eps)}-1\right).$$
Along with $\P_0(T_n^{-\eps}<T_n^{(0,\infty)})=\frac{\tilde W_n(0)}{\tilde W_n(\eps)}$ and $\P_{-\eps}(T_n^{(0,\infty)}<T_n^{-\tau})=1-\frac{\tilde W_n(\eps)}{\tilde W_n(\tau)}$, this gives
\begin{align*}
& \P_0 (\Ups_n-\eps\in\d u,\ T_n^{-\eps}<T_n^{(0,\infty)}<T_n^{-\tau}) \\
=& \frac{\tilde W_n(\tau-\eps)}{\tilde W_n(\tau)}\ \int_{z\in(u,\infty)} N'_n(-\epsilon(\chi-)\in\d u,\ \epsilon(\chi)-\epsilon(\chi-)\in\d z,\ -\inf\epsilon<\tau-\eps,\ \sup\epsilon\geq \eps) \\
=& \frac{\tilde W_n(\tau-\eps)}{\tilde W_n(\tau)}\  \frac n{d_n} e^{-\teta_n u}\d u\ \P_u(T_n^0<T_n^{(\tau-\eps,\infty)}\sachant T_n^0<\infty)\ \int_{z\in(u,\infty)}\tL_n(\d z)  \P_{z-u}(T_n^{(\eps,\infty)}<T_n^0), 
\end{align*}
where in the last equality we first appealed to the strong Markov property at $T_n^{(0,\infty)}$, and then to Proposition \ref{prop_Exc_mesure}. Finally, we get 
\begin{multline*}
 \P_0 (\Ups_n-\eps\in\d u,\ T_n^{-\eps}<T_n^{(0,\infty)}<T_n^{-\tau}) \\
=  \frac n{d_n}\ \frac{\tilde W_n(\tau-u-\eps)}{\tilde W_n(\tau)}\ \d u \ \int_{z\in(u,\infty)}\tL_n(\d z)  \left(1-\frac{\tilde W_n(\eps-(z-u))}{\tilde W_n(\eps)}\right), 
\end{multline*}
which, along with $p'\neps=\frac{d_n}{n}\left(\frac{1}{\tilde W_n(\eps)}-\frac{1}{\tilde W_n(\tau)}\right)$, proves \eqref{formula_nuin}.\\

We next want to prove \eqref{formula_nui}. A similar reasoning as for \eqref{formula_nuin} holds, except that $Z$ has infinite variation. In particular, if $Z$ has a Gaussian component, the process can then creep upwards. Using as before the strong Markov property at $T^{-\eps}$ and Proposition 0.5.2(ii) in \cite{B}, we have 
\begin{align*}
 N'(\Ups&-\eps\in\d u,\ -\inf\epsilon\in[\eps,\tau)) \\
&= N'(-\inf\epsilon>\eps)\P_{-\eps}(T^{(0,\infty)}<T^{-\tau}) \frac{N'(-\epsilon(\chi-)\in\d u,\ -\inf\epsilon<\tau-\eps,\ \sup\epsilon\geq\eps)}{N'(-\inf\epsilon<\tau-\eps,\ \sup\epsilon\geq\eps)}.
\end{align*}
On the one hand, from \cite[Section 4]{Obloj} we know that
\begin{align*}
& N'(-\inf\epsilon>\eps)=\frac{1}{W(\eps)}, \\
& N'(-\inf\epsilon<\tau-\eps,\ \sup\epsilon\geq\eps)=\frac{1}{W(\tau-\eps)}\left(\frac{W(\tau)}{W(\eps)}-1\right),
\end{align*}

which gives
$$\P_{-\eps}(T^{(0,\infty)}<T^{-\tau})\frac{N'(-\inf\epsilon>\eps)}{N'(-\inf\epsilon<\tau-\eps,\ \sup\epsilon\geq\eps)}=\left(1-\frac{W(\eps)}{W(\tau)}\right)\frac{W(\tau-\eps)}{W(\tau)-W(\eps)}=\frac{W(\tau-\eps)}{W(\tau)}.$$

On the other hand, distinguishing the excursions entering $(0,\infty)$ immediately from the others leads to 
\begin{multline*}
 N'(-\epsilon(\chi-)\in\d u,\ -\inf\epsilon<\tau-\eps,\ \sup\epsilon\geq\eps) \\
= N'(-\epsilon(\chi-)\in\d u,\ -\inf\epsilon\in(0,\tau-\eps),\ \sup\epsilon\geq\eps) + N'(-\inf\epsilon=0,\ \sup\epsilon\geq\eps)\delta_0(\d u),
\end{multline*}
where $N'(-\inf\epsilon=0,\ \sup\epsilon\geq\eps)=\frac{b^2}{2}\frac{W'(\eps)}{W(\eps)}$ according to \cite[Section 4]{Obloj}.

Then similarly as before, applying the strong Markov property at $T^{(0,\infty)}$ and Proposition \ref{prop_Exc_mesure}, we get 
\begin{align*}
& N'(\Ups-\eps\in\d u,\ -\inf\epsilon\in[\eps,\tau)) \\
&=\frac{W(\tau-\eps)}{W(\tau)} \Bigg( e^{-\teta u}\d u\ \P_u(T^0<T^{(\tau-\eps,\infty)}\sachant T^0<\infty)\ \int_{z\in(u,\infty)}\L(\d z)  \P_{z-u}(T^{(\eps,\infty)}<T^0) \\
&\qquad\qquad\qquad\qquad\qquad\qquad\qquad\qquad\qquad\qquad\qquad\qquad\qquad\qquad\qquad +\frac{b^2}{2}\frac{W'(\eps)}{W(\eps)}\delta_0(\d u)  \Bigg)\\
&=\frac{W(\tau-\eps)}{W(\tau)} \Bigg( \frac{b^2}{2}\frac{W'(\eps)}{W(\eps)}\delta_0(\d u) + \d u\ \frac{W(\tau-u)}{W(\tau-\eps)}\ \int_{z\in(u,\infty)}\L(\d z)\left(1-\frac{W(\eps-(z-u))}{W(\eps)}\right) \Bigg),
\end{align*}
which proves \eqref{formula_nui}.
\end{demole}

\par\bigskip
\begin{demole}{\ref{lemme_cv_nui}}
We can now prove the weak convergence of $\nuin$ to $\nui$. To begin with, formulas \eqref{formula_nuin} and \eqref{formula_nui} of Lemma \ref{lemme_formula_nui}, along with  the convergence of $\tilde W_n$ towards $W$ (which implies in particular $p'\neps\to p_\eps$ as $n\to\infty$), ensure that we only have to prove the weak convergence of 

$$ \d u\ \mathds1_{u\in(0,\tau-\eps)}\ \tilde W_n(\tau-u-\eps)\ \int_{z\in(u,\infty)}\tL_n(\d z)  \left(1-\frac{\tilde W_n(\eps-(z-u))}{\tilde W_n(\eps)}\right) $$
towards
$$\delta_0(\d u)\frac{b^2}{2}\frac{W'(\eps)}{W(\eps)}W(\tau-\eps) + \d u\ \mathds1_{u\in(0,\tau-\eps)}\ W(\tau-u-\eps)\ \int_{z\in(u,\infty)}\L(\d z)\left(1-\frac{W(\eps-(z-u))}{W(\eps)}\right).$$

First notice that since $\tilde W_n$ and $W$ vanish on the negative half-line, we have 
\begin{align*}
& \int_{(u,\infty)}\tL_n(\d z)  \left(1-\frac{\tilde W_n(\eps-(z-u))}{\tilde W_n(\eps)}\right) =\bar\tL_n(u+\eps)+ \int_{(u,u+\eps]}\tL_n(\d z) \left(1-\frac{\tilde W_n(\eps-(z-u))}{\tilde W_n(\eps)}\right) \\
& \text{and }\int_{(u,\infty)}\L(\d z)\left(1-\frac{W(\eps-(z-u))}{W(\eps)}\right) = \bar\L(u+\eps)+\int_{(u,u+\eps]}\L(\d z)\left(1-\frac{W(\eps-(z-u))}{W(\eps)}\right).
\end{align*}

The functions $u\mapsto \tilde W_n(\tau-u-\eps)\bar\tL_n(u+\eps)$ converge pointwise on $(0,\tau-\eps)$ towards $u\mapsto W(\tau-u-\eps)\bar\L(u+\eps)$ and are bounded by $\sup_{n\geq1} \bar\tL_n(\eps)$. Then by dominated convergence, we have the following weak convergence :
$$\d u\ \tilde W_n(\tau-u-\eps)\bar\tL_n(u+\eps) \mathds1_{u\in(0,\tau-\eps)}\ \Rightarrow\   \d u\  W(\tau-u-\eps) \bar\L(u+\eps)  \mathds1_{u\in(0,\tau-\eps)}.$$
Finally, it remains to prove the weak convergence of
$$\d u \mathds1_{u\in(0,\tau-\eps)} \tilde W_n(\tau-u-\eps) \int_{(u,u+\eps]}\tL_n(\d z) \left(1-\frac{\tilde W_n(\eps-(z-u))}{\tilde W_n(\eps)}\right) $$
towards
$$\delta_0(\d u)\frac{b^2}{2}\frac{W'(\eps)}{W(\eps)}W(\tau-\eps)  + \d u\ \mathds1_{u\in(0,\tau-\eps)}\ W(\tau-u-\eps)\ \int_{z\in(u,u+\eps]}\L(\d z)\left(1-\frac{W(\eps-(z-u))}{W(\eps)}\right).$$

Consider $g$ a continuous bounded function on $\R_+$. We have :
\begin{align*}
& \int_0^{\tau-\eps} \d u\ g(u)\ \tilde W_n(\tau-u-\eps) \left(\int_{(u,u+\eps]}\tL_n(\d z) \left(1-\frac{\tilde W_n(\eps-(z-u))}{\tilde W_n(\eps)}\right)\right) \\
 &= \int_{(0,\tau]} \tL_n(\d z) \left( \int_{0\wedge z-\eps}^z \d u\ g(u) \tilde W_n(\tau-u-\eps)\left(1-\frac{\tilde W_n(\eps-(z-u))}{\tilde W_n(\eps)}\right)\right) \\
 &= \int_{(0,\tau]} \tL_n(\d z) \left( \int_{0}^{z\wedge\eps} \d v\ g(z-v) \tilde W_n(\tau-\eps-z+v)\left(1-\frac{\tilde W_n(\eps-v)}{\tilde W_n(\eps)}\right)\right).
\end{align*}

We set, for all $u\geq0$, 
\begin{center}
$h_n(u):=\int_{0}^{z\wedge\eps} \d v\ g(z-v) \tilde W_n(\tau-\eps-z+v)\big(1-\frac{\tilde W_n(\eps-v)}{\tilde W_n(\eps)}\big)$ \\
$h(u):=\int_0^{z\wedge\eps} \d v\ g(z-v)W(\tau-\eps-z+v)\big(1-\frac{W(\eps-v)}{W(\eps)}\big) $.
\end{center}

We then verify that the conditions of Proposition \ref{prop_Cb_Ku2}, stated in the appendix, are fulfilled :
\begin{itemize}
 \item The functions $h_n$ and $h$ can be bounded by $\eps\cdot\sup|g|$ and are continuous thanks to the continuity on $\R_+$ of the functions $\tilde W_n$ and $W$.
 \item The dominated convergence theorem and the uniform convergence of $\tilde W_n$ towards $W$ on $\R_+$ (see Proposition \ref{prop_cv}.\eqref{prop_cv_W}) ensure that $h_n$ converges uniformly on $\R_+$ towards $h$ (recall that $W(\eps)>0$).
 \item Now since $\tilde W'_n$ converges uniformly towards $W'$ on every compact set of $\R_+^*$ (again from Proposition \ref{prop_cv}.\eqref{prop_cv_W}), the sequence $(v\mapsto\frac1v (\tilde W_n(\eps)-\tilde W_n(\eps-v)))_n$ converges uniformly towards $v\mapsto\frac1v (\tilde W_n(\eps)-\tilde W_n(\eps-v))$ on $(0,\frac\eps2)$. Consequently, for all $a>0$, if $n$ is large enough we have
$$\sup_{u\in(0,\eps/2)}\left|\frac{h_n(u)-h(u)}{u^2}\right|\leq \sup_{u\in(0,\eps/2)}\frac{a}{u^2} \sup|g| \int_0^u v\,\d v =\frac a2 \sup|g|,$$
and thus we have uniform convergence of $u\mapsto h_n(u)/u$ towards $u\mapsto h(u)/u$ on $(0,\eps/2)$.
 \item In the same way we get from the continuity of $g$, $\tilde W_n$ and $W$ that 
$$\frac{h(u)}{u^2}\underset{u\to0}\to \frac12 g(\eps)\frac{W'(\eps)}{W(\eps)}W(\tau-\eps).$$
\end{itemize}
We then get the expected convergence from an appeal to Proposition \ref{prop_Cb_Ku2}. As a conclusion, we proved that the measures $\nuin$ converge weakly to $\nui$.
\end{demole}

\par\bigskip

\begin{lemme} \label{lemme_cv_M0}
 As $n\to\infty$, $M\neps(0)$ converges in distribution towards $M_\eps(0)$.
\end{lemme}

\begin{demo}
 We have to prove the weak convergence of 
\begin{center}
$\nuin(\d u\times\un)$, \quad$\int_{[\eps,\tau)} \nuin(\d x\times\zero) \num_n(x,\d u)$\quad and \quad$\int_{[\eps,\tau)} \nuin(\d x,\zero)\nud_n(x,\d u)$.                                                                        
\end{center} 
The convergence of the first one is a straightforward consequence of Lemma \ref{lemme_cv_nui}. We prove below the convergence of the Laplace transform of the second one, and a similar reasoning holds for the third one. We consider for all $a\geq0$
\begin{align*}
\int_{(0,\infty)} e^{-au} \int_{[\eps,\tau)} \nuin(\d x\times\zero) \num_n(x,\d u)=
 \int_{[\eps,\tau)} \nuin(\d x\times\zero) \int_{[\eps,\tau)} e^{-au}\num_n(x,\d u) 
\end{align*}
and set $h_n(x):=\int_{[\eps,\tau)} e^{-au}\num_n(x,\d u) $ and $h(x):=\int_{[\eps,\tau)} e^{-au}\num(x,\d u)$.\\

The functions $h$ and $h_n$ are all bounded by $1$. Moreover, they are continuous : indeed, we have 
$$|\E(e^{-a\Ha(\e)},\ \Linv(\e)<\Tx)-\E(e^{-a\Ha(\e)},\ \Linv(\e)<T^{-x_0})|\leq\P(T^{-x_0}<\Tx),$$
which vanishes as $x\to x_0$ thanks to the a.s. continuity of $x\mapsto\Tx$ on $\R_+$ under $\P$. Here again, for the sake of simplicity, we omitted the conditioning, but a similar reasoning and an appeal to the continuity of $W$, lead to the continuity of $h$. Besides, the same arguments can be used to get the continuity of $h_n$. Finally, as established in the proof of \cite[Th. 4]{Karr}, Lemma \ref{lemme_cv4b_Karr} ensures the uniform convergence of $h_n$ towards $h$ on every compact set of $\R_+^*$. Then, since $\nuin$ and $\nui$ are probability measures such that $\nuin\Rightarrow\nui$, Lemma \ref{lemme_Cb} entails the convergence of the Laplace transform of $\int_{[\eps,\tau)} \nuin(\d x\times\zero) \num_n(x,\d u)$ towards that of $\int_{[\eps,\tau)} \nui(\d x\times\zero) \num(x,\d u)$.
\end{demo}

\par\bigskip
\begin{demopr}{\ref{prop_cv_M}}
 Lemma \ref{lemme_cv4b_Karr} ensures that the Markov chains $M\neps$ and $M_\eps$ satisfy condition (4).b in \cite{Karr}, while condition (4).a of the same paper is given by Lemma \ref{lemme_cv_M0}. Then the announced convergence is a consequence of \cite[Theorem (1)]{Karr}.
\end{demopr}

\section{Appendix}

\subsection{A convergence lemma for integrals}

\begin{lemme}\label{lemme_Cb}
 Let $(h_n)_{n\geq0}$ and $h$ be continuous bounded mappings from $\R^d$ to $\R$, and suppose $(h_n)$ is dominated by a bounded function. Let $(\mu_n)_{n\geq0}$ and $\mu$ be in $\mathcal M_f(\R^d)$ and suppose that
\begin{enumerate}[\upshape(i)]
 \item $(\mu_n)$ converges weakly to $\mu$.
 \item The sequence of mappings $(h_n)$ converges to $h$ uniformly on every compact set of $\R^d$.
\end{enumerate}
Then \upshape $$\int h_n\d\mu_n \xrightarrow{n\to\infty} \int h\d\mu$$
\end{lemme}

\begin{demo}
We have
$$  \left|\int h_n\d\mu_n-\int h\d\mu\right|\leq \left|\int (h_n-h)\d(\mu_n-\mu)\right|+ \left|\int (h_n-h)\d\mu\right|+ \left|\int h\d(\mu_n-\mu)\right|.$$
The mapping $h$ is continuous and bounded on $\R^d$, then (i) implies the convergence to $0$ of the term $|\int h\d(\mu_n-\mu)|$. The domination and convergence assumptions made on $(h_n)$ allow us to apply the dominated convergence theorem to get the convergence of $|\int (h_n-h)\d\mu|$ to $0$. As for the first term in the sum, it requires some additional details : Let $\eps$ be a positive real number. First, thanks to (i) and since $(h_n-h)$ is dominated by a constant, we can find a compact set $K_\eps\subset \R^d$ and $n_0\in\N$ such that $|\int_{K_\eps^c}(h_n-h)\d(\mu_n-\mu)|\leq\eps$ for $n\geq n_0$. Secondly the uniform convergence on the compact set $K_\eps$ of the sequence $(h_n)$ ensures that $|\int_{K_\eps}(h_n-h)\d(\mu_n-\mu)|\leq\eps$ for $n$ large enough. In consequence we have convergence of the term $|\int (h_n-h)\d(\mu_n-\mu)|$ to 0, and the result follows.
\end{demo}

\subsection{Consequences of Assumption \HypLevy}

Then we turn our attention to the proof of Proposition \ref{prop_cv}. \par\medskip

\noindent\emph{Proof of Proposition \ref{prop_cv} :} 
\begin{enumerate}[(i)]
 \item Since $Z$ is a.s. continuous at $\Tx$ (resp. is not a compound Poisson process), we have $\underset{\eps\to0+}\lim T^{-(x+\eps)}=\Tx$ (resp.$\underset{\eps\to0+}\lim T^{(y+\eps,\infty)}=T^{(y,\infty)}$) a.s., and hence the convergence in law of $\Txn$ towards $\Tx$ (resp. $T_n^{(y,\infty)}$ towards $T^{(y,\infty)}$) is a straightforward consequence of Proposition VI.2.11 in \cite{JS}.
 
\item Now $\phi_n$ (resp. $\phi$) is the Laplace exponent of the process $x\mapsto \Txn$ (resp. $x\mapsto \Tx$) \cite[Th. VII.1.1]{B}. The pointwise convergence of $\tilde\phi_n$ to $\phi$ is thus a consequence of point (\ref{prop_cv_Tx}). The uniform convergence comes from the fact that for all $n\geq1$, $\tilde\phi_n$ is increasing on $\R_+$.
 
\item The proof of the pointwise convergence of $\tilde W_n$ (resp. $\tilde W'_n$) towards $W$ (resp. $W'$) can be found in \cite[Prop. 3.1]{LS}. Moreover, we have for all $y>x$ $\P(T^{-x}<T^{(y-x,\infty)})=\frac{W(x)}{W(y)}$ \cite[Th. VII.2.8]{B}, and then the function $x\mapsto \tilde W_n(x)/\tilde W_n(y)$ is decreasing, thus the convergence of $\tilde W_n$ towards $W$ is uniform on every compact set of $\R_+$. Finally, the uniform convergence of $\tilde W'_n$ towards $W'$ on every compact set of $\R_+^*$ can be deduced from the expression of $\tilde W'_n$ given in the proof of Lemma 8.2 in \cite{K}, as a product of two monotone functions.
\end{enumerate}\hfill $\square$

Finally, we recall the following result, obtained in \cite{LPWM1} as a consequence of Assumption \HypLevy.

\begin{proposition}\label{prop_Cb_Ku2}
 Let $(g_n)_{n\geq0}$ and $g$ be continuous bounded mappings from $\R_+$ to $\R$, where $g$ satisfies $g(u)/u^2 \to K$ as $u\to 0+$ for some constant $K$. Assume that the mappings $\tilde g_n:u\mapsto \frac{g_n(u)}{1\wedge u^2}$ converge uniformly to  $\tilde g:u\mapsto \frac{g(u)}{1\wedge u^2}$ on $\R_+^*$. Then as $n\to\infty$,
\upshape $$\tL_n(g_n) \underset{n\to\infty}{\to} \L(g) + Kb^2.$$
\end{proposition} 

\begin{remerciements}
 I would like to thank my supervisor, Amaury Lambert, for his very helpful advice and encouragement.
\end{remerciements}

\selectlanguage{english}
\par\bigskip
\nocite{*}
\bibliographystyle{plain}
\bibliography{biblioLPWM2}

\end{document}